\documentclass[12pt]{amsart}
\usepackage{mathrsfs}
\usepackage{amsfonts}
\usepackage{eufrak}
\usepackage{amsmath}
\usepackage{upgreek}
\usepackage{amssymb}
\usepackage{bm}
\usepackage{latexsym}
\usepackage{verbatim}
\usepackage[margin=3.5cm]{geometry}

\newcommand{\bdelta}{\boldsymbol\delta}
\newcommand{\bgamma}{\boldsymbol\gamma}
\newcommand{\done}{\hfill$\square$}
\newcommand{\ord}{\mbox{ord}}
\newcommand{\PGL}{\mbox{PGL}}

\begin{document}

\title{Eisenstein Series in Kohnen Plus Space for Hilbert Modular Forms}
\author{Ren He Su}
\address{Graduate school of mathematics, Kyoto University, Kitashirakawa, Kyoto, 606-8502, Japan}
\email{ru-su@math.kyoto-u.ac.jp}
\thanks{The author thanks Prof. ~Ikeda for his advice.}
\maketitle
\begin{abstract}

In 1975, Cohen constructed a kind of one-variable modular forms of half-integral weight, says $r+(1/2),$ whose $n$-th Fourier coefficient $H(n)$ only occurs when $(-1)^r n$ is congruent to 0 or 1 modulo 4.
The space of modular forms whose Fourier coefficients have the above property is called Kohnen plus space, initially introduced by Kohnen in 1980.
Recently, Hiraga and Ikeda generalized the Kohnen plus space to the spaces for half-integral weight Hilbert modular forms with respect to general totally real number fields.
If one such Hilbert modular form $f$ of parallel weight $\kappa+(1/2)$ lying in a generalized Kohnen plus space has $\xi$-th Fourier coefficients $c(\xi)$, then $c(\xi)$ does not vanish only if $(-1)^\kappa\xi$ is congruent to a square modulo 4.
In this paper, we use an adelic way to construct Eisenstein series of parallel half-integral weight belonging to the generalized Kohnen plus spaces and give an explicit form for their Fourier coefficients.
These Eisenstein series give a generalization of the modular forms introduced by Cohen.
Moreover, we show that the Kohnen plus space is generated by the cusp forms and the Eisenstein series we constructed as a vector space over $\mathbb{C}.$

\end{abstract}


\section*{Introductions}
In the beginning of this paper, the definitions and some properties of the generalized Kohnen plus space for Hilbert modular forms of half-integral weights will be introduced. Then we give the idea and the main results of this paper.\\
\hbox{}\quad Let $F$ be a totally real number field over $\mathbb{Q}$ with degree $n$ while $\mathfrak{o}, \mathfrak{d}$ are its ring of integers and the different over $\mathbb{Q}.$
Denote the $n$ embeddings of $F$ to $\mathbb{R}$ by $\iota_1, \iota_2, ..., \iota_n.$
The congruence subgroup $\Gamma$ of $SL_2(\mathfrak{F})$ is defined by
\[
\Gamma=
\left\{
\begin{pmatrix}a&b\\c&d\end{pmatrix}
\in SL_2(F)\,\bigg|\,a, d\in\mathfrak{o}, c\in\mathfrak{d}^{-1}, d\in4\mathfrak{d}
\right\}.
\]
Now if we set the theta function $\theta$ on $\mathfrak{h}^n$ by
\[
\theta(z)=\sum_{\xi\in\mathfrak{o}}q^\xi
\mbox{ for }z=(z_1,...,z_n)\in\mathfrak{h}^n
\]
where $q^\xi=\exp(2\pi\sqrt{-1}\sum_{i=1}^n\iota_i(\xi)z_i),$ then for any $\gamma\in\Gamma$ we have
\[
\theta(\gamma z)=j^{1/2}(\gamma,z)\theta(z)
\]
for some factor of automorphy $j^{1/2}(\gamma,z).$
Here $\gamma z$ is the usual M\"{o}bius transformation.\\
\hbox{}\quad Let $\kappa\geq1$ be some integer.
Write $j^{\kappa+1/2}(\gamma,z)=(j^{1/2}(\gamma,z))^{2\kappa+1}.$ Let $M_{\kappa+1/2}(\Gamma)$ ($S_{\kappa+1/2}(\Gamma)$) be the space of consisting of Hilbert modular forms (cusp forms) on $\Gamma$ of parallel weight $\kappa+1/2$ with respect to the factor of automorphy $j^{\kappa+1/2}(\gamma,z).$\\
\hbox{}\quad Now set $F=\mathbb{Q}$.
For any nonzero integer $n,$ denote $\chi_n$ the character of the quadratic extension $\mathbb{Q}(\sqrt{n})/\mathbb{Q},$ $D_n$ the discriminant of $\mathbb{Q}(\sqrt{n})/\mathbb{Q}$ and $f_n$ the positive integer such that $n=D_nf_n^2.$
Then for $r\geq2,$ Cohen \cite{Cohen:75} introduced the special modular form $\mathcal{H}_r\in M_{r+1/2}(\Gamma)$ of weight $r+1/2$ which has the Fourier expansion
\[
\begin{aligned}
&\mathcal{H}_r(z)=\zeta(1-2r)+\\
&\sum_{\begin{smallmatrix}N\geq1\\(-1)^rN\equiv0,1\mathrm{ mod }4\end{smallmatrix}}\left(L(1-r,\chi_{(-1)^rN})\sum_{d\,|\,f_{(-1)^rN}}\mu(d)\chi_{(-1)^rN}(d)d^{r-1}\sigma(f_{(-1)^rN}/d)\right)q^N.
\end{aligned}
\]
\hbox{}\quad Back to general $F.$
For any $\xi\in F,$ we denote $\xi\equiv\square$ mod $4$ if there exists $x\in\mathfrak{o}$ such that $\xi-x^2\in4\mathfrak{o}.$
Notice that each $h\in M_{\kappa+1/2}(\Gamma)$ can be written in a Fourier expansion as
\[
h(z)=\sum_{\xi\in\mathfrak{o}}c(\xi)q^\xi.
\]
The Kohnen plus space $M^+_{\kappa+1/2}(\Gamma)$ with respect to $M_{\kappa+1/2}(\Gamma)$ is the subspace consisting of $h$ whose $\xi$-th Fourier coefficient disappears unless $(-1)^\kappa\xi\equiv\square$ mod $4$. That is,
\[
M^+_{\kappa+1/2}(\Gamma)=\left\{
h\in M_{\kappa+1/2}(\Gamma)\,\bigg|\,h(z)=\sum_{(-1)^\kappa\xi\equiv\square\mbox{ mod }4}c(\xi)q^\xi\right\}
\]
Also, we set $S^+_{\kappa+1/2}(\Gamma)=M^+_{\kappa+1/2}(\Gamma)\cap S_{\kappa+1/2}(\Gamma).$ It is called the Kohnen plus space, too.\\
\hbox{}\quad Let us first consider the case $F=\mathbb{Q}.$
The Kohnen plus space under this condition was first introduced by Kohnen \cite{Kohnen:80} in 1980.
Obviously the modular form $\mathcal{H}_r$ mentioned above is in it.
If for any $h\in M_{\kappa+1.2}(\Gamma)$ we set
\[
\begin{aligned}
(U_4h)(z)&=\frac{1}{4}\sum_{j=1}^{4}h\left(\frac{z+j}{4}\right),\\
(W_4h)(z)&=(-2iz)^{-\kappa-1/2}f(-\frac{1}{4z}),
\end{aligned}
\]
then the two operators leave $S_{\kappa+1/2}(\Gamma)$ stable. Moreover, Kohnen showed that $S^+_{\kappa+1/2}(\Gamma)$ is the eigenspace for the operator $W_4U_4$ with eigenvalue $(-1)^{\kappa(\kappa+1)/2}2^\kappa.$\\
\hbox{}\quad Kohnen also showed that $S^+_{\kappa+1/2}(\Gamma)$ and $M^+_{\kappa+1/2}(\Gamma)$ are isomorphic to $S_{2\kappa}(SL_2(\mathbb{Z}))$ and $M_{2\kappa}(SL_2(\mathbb{Z}))$ as Hecke algebras, respectively. Moreover, if for any $h(z)=\sum_{n\geq0}c(n)q^n$ in $M^+_{\kappa+1/2}(\Gamma)$ we set
\[
\begin{aligned}
&(T^+_{\kappa+1/2}(p^2)h)(z)=\sum_{\begin{smallmatrix}n\geq0\\(-1)^\kappa n\equiv0,1\mathrm{ mod }4\end{smallmatrix}}
\left(
c(p^2n)+\left(\frac{(-1)^\kappa n}{p}\right)p^{\kappa-1}c(n)+p^{2\kappa-1}c\left(\frac{n}{p^2}\right)
\right)q^n\\
&(c(x)=0\mbox{ if }x\notin\mathbb{Z}_{\geq0})
\end{aligned}
\]
where $p$ is any prime number, then $S^+_{\kappa+1/2}(\Gamma)$ has an orthogonal basis of common eigenforms for all $T^+_{\kappa+1/2}(p^2).$\\
\hbox{}\quad For general totally real number field $F$ the Kohnen plus space was introduced by Hiraga and Ikeda \cite{IkedaHiraga} in 2013. Some analogues of the properties for the case $F=\mathbb{Q}$ stated above are showed. For any non-archimedean place $v$ of $F,$ we let $\widetilde{SL_2(F_v)}$ be the metaplectic double covering of $SL_2(F_v).$ Set
\[
\Gamma_v=
\left\{
\begin{pmatrix}a&b\\c&d\end{pmatrix}
\in SL_2(F_v)\,\bigg|\,a, d\in\mathfrak{o_v}, c\in\mathfrak{d_v}^{-1}, d\in4\mathfrak{d_v}
\right\}
\]
and $\widetilde{\Gamma_v}$ be its inverse image in $\widetilde{SL_2(F_v)}.$
There exists a genuine character $\varepsilon_v:\widetilde{\Gamma_v}\rightarrow\mathbb{C}^\times$ defined from the Weil representation of $\widetilde{SL_2(F_v)}.$
Denote $\widetilde{\mathcal{H}}_v=\widetilde{\mathcal{H}}(\widetilde{\Gamma_v}\backslash \widetilde{SL_2(F_v)}/\widetilde{\Gamma_v},\varepsilon_v)$ the Hecke algebra consisting of all compactly supported function $f$ on $\widetilde{SL_2(F_v)}$ such that $f(h_1kh_2)=\varepsilon_v(h_1)\varepsilon_v(h_2)f(k)$ for any $h_1, h_2\in\widetilde{\Gamma_v}.$
Hiraga and Ikeda defined some idempotent $E^K_v\in\widetilde{\mathcal{H}}_v$ such that $E^K=\prod_{v<\infty}E^K_v$ is the projection of $M_{\kappa+1/2}(\Gamma)$ into $M^+_{\kappa+1/2}(\Gamma).$ In particular, if $F=\mathbb{Q},$ since the only other eigenvalue of $W_4U_4$ is $-(-1)^{\kappa(\kappa+1)/2}2^{\kappa-1},$ we have
\[
E^K=(-1)^{\kappa(\kappa+1)/2}\frac{2^{1-\kappa}}{3}W_4U_4+\frac{1}{3}.
\]
\hbox{}\quad Denote $\mathbb{A}$ the adele ring of $F$ and $\mathcal{K}_0=\prod_{v<\infty}PGL_2(\mathfrak{o}_v)$ an open compact subgroup of $PSL_2(\mathbb{A}).$ Let $A^\mathrm{cusp}_{2\kappa}(PGL_2(F)\backslash PSL_2(\mathbb{A})/\mathcal{K}_0)$ be the space of automorphic forms on $PGL_2(F)\backslash\widetilde{SL_2(\mathbb{A})}$ of weight $2\kappa,$ which is invariant by $\mathcal{K}_0.$
Then Hiraga and Ikeda showed that $S^+_{\kappa+1/2}(\Gamma)$ and $A^\mathrm{cusp}_{2\kappa}(PGL_2(F)\backslash PSL_2(\mathbb{A})/\mathcal{K}_0)$ have the same dimension as vector spaces over $\mathbb{C}$ if $\kappa\geq2.$
Moreover, $S^+_{\kappa+1/2}(\Gamma)$ is a direct sum of modular forms which are eigenforms with respect to $\widetilde{H}_v$ for all finite places $v$ not dividing $2$.\\
\hbox{}\quad The main aim of this paper is to generalize the modular form $\mathcal{H}_r$ mentioned above to the cases for general Hilbert modular forms, that is, to find an element in general Kohnen plus space which corresponds to $\mathcal{H}_r.$ Recall that $\mathcal{H}_r$ is only for $r\geq2.$ In fact, if $F\neq\mathbb{Q},$ we can find a Hilbert modular form of weight $3/2$ which has a similar form of Fourier coefficients, while in the case $F=\mathbb{Q}$ there only exists a non-analytic function transforming under $\Gamma$ like a modular form of weight $3/2.$
We call the generalization Eisenstein series in the Kohnen plus space.
The following theorem will be showed in Section 9 (Theorem 9.1) and Section 10 (Theorem 10.3).\\[0.2cm]
\textbf{Theorem.}
We use the same notations as above and let $\kappa$ be an positive integer which is not equal to $1$ if $F=\mathbb{Q}.$ Also, we set $\chi'$ a character of the class group of $F.$ Then we have $G(z)=G_{\kappa+1/2}(z,\chi')\in M^+_{\kappa+1/2}(\Gamma)$ which is defined by
\[
G(z)=L_F(1-2\kappa,\overline{\chi'}^2)+\sum_{\begin{smallmatrix}(-1)^{\kappa}\xi\equiv\square\,\mathrm{mod}\,4\\ \xi\succ 0\end{smallmatrix}}\chi'(\mathfrak{D}_{(-1)^\kappa\xi})L_F(1-\kappa,\overline{\chi_{(-1)^\kappa\xi}\chi'})\mathfrak{C}_\kappa((-1)^\kappa\xi)q^\xi.
\]
where
\[
\mathfrak{C}_\kappa(\xi)
=\sum_{\mathfrak{a}|\mathfrak{F}_\xi}\mu(\mathfrak{a})\chi_\xi(\mathfrak{a})\chi'(\mathfrak{a})N_{F/\mathbb{Q}}(\mathfrak{a})^{\kappa-1}\sigma_{2\kappa-1,\chi'^2}(\mathfrak{F}_\xi\mathfrak{a}^{-1}).
\]
Here $\mathfrak{D}_\xi$ is the relative discriminant of $F(\sqrt{\xi})/F,$ $\mathfrak{F}_\xi^2\mathfrak{D_\xi}=(\xi),$ $\mathfrak{a}$ runs over all integral ideals dividing $\mathfrak{F_\xi},$ $\mu$ is the M\"{o}bius function for ideals and $\sigma_{k,\chi}$ is the sum of divisors function twisted by $\chi,$ that is,
\[
\sigma_{k,\chi}(\mathfrak{A})=\sum_{\mathfrak{b}|\mathfrak{A}}N_{F/\mathbb{Q}}(\mathfrak{b})^k\chi(\mathfrak{b})
\]
for any integral ideal $\mathfrak{A}$ of $F.$
Moreover, $G$ is a Hecke eigenform with respect to the Hecke algebra $\widetilde{\mathcal{H}_v}$ for any finite $v$ which is not even.\\[0.2cm]
\hbox{}\quad So for any $F$ with class number $h$ we can get $h$ Eisenstein series in $M^+_{\kappa+1/2}(\Gamma).$ Using this fact we will show the following corollary as the second main result of this paper in Section 11 (Corollary 11.3).\\[0.2cm]
\textbf{Corollary.}
The Kohnen plus space $M^+_{\kappa+1/2}(\Gamma)$ is a vector space over $\mathbb{C}$ spanned by cusp forms and the $h$ Eisenstein series we got in the last theorem, that is,
\[
M^+_{\kappa+1/2}(\Gamma)=S^+_{\kappa+1/2}(\Gamma)\oplus\bigoplus_{j=1}^{h}\mathbb{C}\cdot G_{\kappa+1/2}(z,\chi_j)
\]
where $\chi_1, ...,\chi_h$ are the $h$ distinct characters of the class group of $F.$\\[0.2cm]
\hbox{}\quad Together with the results of Ikeda and Hiraga we get that $M^+_{\kappa+1/2}(\Gamma)$ is a direct sum of spaces spanned by eigenforms with respect to the Hecke algebra $\widetilde{\mathcal{H}_v}$ for any finite $v$ which is not even.\\
\hbox{}\quad Now let us show the content of this paper. In Section 1, we recall some basic properties of the metaplectic group and Weil representation on it. In Section 2, we define two invariants which will be used later. In Section 3, we introduce the idempotent Hecke operators $e^K$ and $E^K$ and calculate some local integrals. In Section 4, we define the function $f^+_K$ and calculate some values of its Whittaker function. $f^+_K$ will be used in the construction of the Eisenstein series. In Section 5, we consider the Archimedean case. In Section 6, we define the metaplectic group $\widetilde{SL_2(\mathbb{A})}$ and define the automorphic forms on it. We connect the automorphic forms and the Hilbert modular forms in $M_{\kappa+1/2}(\Gamma).$ In Section 7, we define the Kohnen plus space for Hilbert modular forms. In Section 8, we construct the Eisenstein series. In Section 9 and Section 10, we calculate the Fourier coefficients of the Eisenstein series for the case $\kappa\geq2$ and $\kappa=1,$ respectively. In Section 11, we prove the corollary stated above. Finally, in Section 12, we give an example.\\

\section*{Notations}
As usual, the upper-half plane $\mathfrak{h}$ is defined as $\mathfrak{h}=\{z\in\mathbb{C}\,|\,\mbox{Im}(z)>0\}.$\\
\hbox{}\quad For $\alpha\in\mathbb{C},$ we set $\textbf{e}(\alpha)=\exp(2\pi\sqrt{-1}\alpha).$
For a totally real number field $F$ with degree $n$ , if $\iota_1, \iota_2, ..., \iota_n$ are the $n$ embeddings of $F$ into $\mathbb{R},$ we put $\mathbf{e}(\xi\ast z)=\exp(2\pi\sqrt{-1}\sum_{i=1}^{n}\iota_i(\xi)z_i)=\prod_{i=1}^n\textbf{e}(\iota_i(\xi)z_i),$
where $\xi\in F$ and $z=(z_1,z_2,...,z_n)\in\mathfrak{h}^n.$
In particular if the letter $z$ is used, we tend to write $q^\xi=\mathbf{e}(\xi\ast z)$.\\
\hbox{}\quad If $\alpha,\,\beta\in\mathbb{C}$ and $\alpha\neq 0,$ we define $\alpha^\beta=\exp(\beta\log(\alpha)),$ where $-\pi<\mbox{Im}(\log(\alpha))\leq\pi.$\\


\section{The Metaplectic Cover and The Weil Representation}
Let $F$($\neq\mathbb{C}$) be a local field with characteristic $0$.
If $F$ is non-archimedean and a finite extension over $\mathbb{Q}_p$, $\mathfrak{o},$ $\mathfrak{p}$ and $\mathfrak{d}$ denote the ring of integers, the prime ideal and the different over $\mathbb{Q}_p.$ 
Also, let $q$ be the order of the residue field and $\varpi$ is a prime element.
On the other hand, if $F$ is real then $\mathfrak{d}$ also denotes the different over $\mathbb{Q}.$
Furthermore, we set $\psi :F\rightarrow \mathbb{C}^\times$ to be a non-trivial additive character on $F$ such that $\psi(x)=\mathbf{e}(x)$ if $F\simeq\mathbb{R}.$
If F is non-archimedean, let $c_\psi$ denote the order of $\psi,$ that is, $c_\psi$ is the largest non-negative integer $c$ such that $\psi(\mathfrak{p}^{-c})=1.$
In this paper, we take $\psi$ so that $c_\psi$ is equal to the order of $\mathfrak{d}.$
We let $\boldsymbol\delta=1$ if $F\simeq\mathbb{R}$, and an element with order $c_\psi$ if $F$ is non-archimedean.
The Haar measure $dx$ is the usual Lebesgue measure if $F\simeq\mathbb{R}$, and normalized so that Vol($\mathfrak{o}$)$=1$ if $F$ is non-archimedean.\\
\hbox{}\quad
We denote $\widetilde{{SL}_2(F)}$ the metaplectic covering of group of $SL_2(F),$ that is,
\[\widetilde{SL_2(F)}=\{[g,\zeta]|g\in SL_2(F),\zeta\in\{\pm 1\}\}\]
equipped with the multiplication
\[[g_1,\zeta_1]\cdot[g_2,\zeta_2]=[g_1g_2,\zeta_1\zeta_2\mathbf{c}(g_1,g_2)],\]
where
\[
\mathbf{c}(g_1,g_2)=\langle\frac{\tau(g_1)}{\tau(g_1g_2)},\frac{\tau(g_2)}{\tau(g_1g_2)}\rangle,
\qquad
\tau\left(\begin{pmatrix}a&b\\c&d\end{pmatrix}\right)=
\begin{cases}
c&\mbox{ if }c\neq0,\\
d&\mbox{ if }c=0.
\end{cases}
\]
Here $\langle\;,\;\rangle$ denotes the Hilbert symbol and $\mathbf{c}$ is called the Kubota's 2-cocycle on $SL_2(F)$. For any subset $S$ of $SL_2(F)$, $\widetilde{S}$ denotes its inverse image in $\widetilde{SL_2(F)}.$ Now for any $g\in SL_2(F),$ we denote $[g,1]$ by $[g]$ and set
\[ 
\mathbf{u}^\sharp(x)=\left[\begin{pmatrix}1&x\\0&1\end{pmatrix}\right],\qquad \mathbf{u}^\flat(x)=\left[\begin{pmatrix}1&0\\x&1\end{pmatrix}\right],\qquad\mbox{for } x\in F,
\]
\[
 \mathbf{m}(a)=\left[\begin{pmatrix}a&0\\0&a^{-1}\end{pmatrix}\right],\qquad \mathbf{w}_a=\left[\begin{pmatrix}0&-a^{-1}\\a&0\end{pmatrix}\right],\qquad\mbox{for } a\in F^\times .
\]
\quad
For each Schwartz function $\phi\in\mathcal{S}(F)$, its Fourier transform is defined by
\[\hat{\phi}(x)=|\boldsymbol\delta|^{1/2}\int_F \phi(y)\psi(xy)dy.\] 
\quad
Here $|\boldsymbol\delta|^{1/2}dy$ is the self-dual Haar measure for the Fourier transform.
Under this definition, for any $\phi\in\mathcal{S}(F)$ and $a\in F^\times,$ we have
\begin{equation}
\label{wc}
\int_F\phi(x)\psi(ax^2)dx=\alpha_\psi(a)|2a|^{-1/2}\int_F\hat{\phi}(x)\psi\left(-\frac{x^2}{4a}\right)dx
\end{equation}
for the so-called Weil constant $\alpha_\psi(a).$ The Weil constant has the properties that $\alpha_\psi(a)^8=1,$ $\alpha_\psi(b^2a)=\alpha_\psi(a)$ for all $b\in F^\times,$ $\alpha_\psi(-a)=\overline{\alpha_\psi(a)}.$
Also, for any $a,\,b\in F^\times,$ we have
\[\frac{\alpha_\psi(a)\alpha_\psi(b)}{\alpha_\psi(1)\alpha_\psi(ab)}=\langle a,b\rangle.\]
\quad
The Weil representation $\omega_\psi$ of $\widetilde{SL_2(F)}$ on $\mathcal{S}(F)$ is defined by
\[
\omega_\psi(\mathbf{m}(a))\phi(t)=\frac{\alpha_\psi(1)}{\alpha_\psi(a)}|a|^{1/2}\phi(at),
\]
\[
\omega_\psi(\mathbf{u}^\sharp(b))\phi(t)=\psi(bt^2)\phi(t),
\]
\[
\omega_\psi(\mathbf{w}_a)\phi(t)=\overline{\alpha_\psi(a)}|2a^{-1}|^{1/2}\hat{\phi}(-2a^{-1}t).
\]
Noticing $\mathbf{u}^\flat(c)=\mathbf{u}^\sharp(c^{-1})\mathbf{w}_c\mathbf{u}^\sharp(c^{-1}),$ one easily deduce that for $c\in F^\times,$
\[
\omega_\psi(\mathbf{u}^\flat(c))\phi(t)=\overline{\alpha_\psi(c)}|2\boldsymbol\delta c^{-1}|^{1/2}\int_F \phi(t+y)\psi(c^{-1}y^2)dy.
\]
The Weil representation $\omega_\psi$ is unitary with the inner product
\[
(\phi_1,\phi_2)=\int_F \phi_1(t)\overline{\phi_2(t)}dt
\]
for any $\phi_1,\,\phi_2\in\mathcal{S}(F).$\\
\hbox{}\quad
For any fractional ideals $\mathfrak{a}$ and $\mathfrak{b}$ of $F$ such that $\mathfrak{ab}\subset \mathfrak{o},$ we define a compact open subgroup $\Gamma[\mathfrak{a,b}]<SL_2(F)$ by
\[
\Gamma[\mathfrak{a,b}]
=\left\{
\begin{pmatrix}a&b\\c&d\end{pmatrix}\bigg|
a,\,b\in \mathfrak{o},\,b\in\mathfrak{a},\,c\in\mathfrak{b}
\right\}
\]
and set
\[
\Gamma=\Gamma[\mathfrak{d}^{-1},4\mathfrak{d}].
\]
\textbf{Lemma 1.1.}
Let $\phi_0\in\mathcal{S}(F)$ be the characteristic function of $\mathfrak{o},$ then Lemma 1.1 of \cite{IkedaHiraga} shows that there is a genuine character $\varepsilon:\widetilde{\Gamma}\rightarrow\mathbb{C}^\times$ such that
\[
\omega_\psi(g)\phi_0=\varepsilon^{-1}(g)\phi_0
\]
for any $g\in\widetilde{\Gamma}.$
If $q$ is even, then
\[
\varepsilon\left(\left[
\begin{pmatrix}a&b\\c&d\end{pmatrix},\zeta
\right]\right)=
\begin{cases}
\zeta\alpha_\psi(d)\overline{\alpha_\psi(1)}&\mbox{if }c=0,\\
\zeta\alpha_\psi(1)\overline{\alpha_\psi(d)}\langle c,d\rangle &\mbox{if }c\neq0.
\end{cases}
\]
\hbox{}\quad In addition, it is worth mentioning that
\[
\hat{\phi_0}(x)=|\boldsymbol\delta|^{1/2}\phi_0(\boldsymbol\delta x).
\]


\section{Definitions of $\mathfrak{f}_\xi$ and $\chi_\xi$}
From now on we let $F$ be a non-archimedean local field and use the same notations as in the last section until Section 4. Let $e$ be the order of $2$ in $F.$\\[0.2cm]
\textbf{Definition 2.1.} For a non-negative integer $r,$ set
\[
U_r=
\left\{\begin{array}{lr}
\mathfrak{o}^\times\\
1+\mathfrak{p}^r\mathfrak{o}
\end{array}\right.
\begin{array}{lr}
\mbox{if }r=0,\\
\mbox{if }r>0,
\end{array}
\]
\[
\mathcal{U}_r=
\left\{\begin{array}{lr}
\mathfrak{o}^{\times 2}U_{2r}\\
\emptyset&
\end{array}\right.
\begin{array}{lr}
\mbox{if }0\leq r\leq e,\\
\mbox{if }r>e.
\end{array}
\]
\quad
We have the following two well-known lemmas:\\[0.2cm]
\textbf{Lemma 2.2.} The following assertion holds.
\[
\mathfrak{o}^{\times 2}U_{2r}=\mathfrak{o}^{\times 2}U_{2r+1}=\mathcal{U}_r \mbox{ for } 0\leq r<e.
\]
\textbf{Lemma 2.3.} We denote $d_\xi$ the order of the conductor of $F(\sqrt{\xi})/F$ for $\xi\in F^\times.$ If $\xi\in\mathcal{U}_r\setminus\mathcal{U}_{r+1}$ ($0\leq r\leq e$), then $d_\xi=2e-2r.$ If $\ord(\xi)$ is odd, then $d_\xi=2+1.$\\[0.2cm]
\hbox{}For the proofs, see \cite{O'Meara:63} \textsection \, 63A and \cite{Okazaki:91} Proposition 3.
\\ \hbox{}\quad
Now we have our definitions of $\mathfrak{f}_\xi$ and $\chi_\xi$:\\[0.2cm]
\textbf{Definition 2.4.} For $\xi\in F^{\times},$ we define
\[
\begin{array}{lr}
\mathfrak{f}_\xi=(\ord(\xi)-d_\xi)/2,
\\[0.2cm]
\chi_\xi=
\left\{\begin{array}{lr}
1 \\ -1 \\ 0
\end{array}\right.
\begin{array}{lr}
\mbox{if }\xi\in F^{\times 2},\\
\mbox{if }F(\sqrt{\xi})/F\mbox{ is an unramified quadratic extension,}\\
\mbox{if }F(\sqrt{\xi})/F\mbox{ is a ramified quadratic extension.}
\end{array}
\end{array}
\]
We also set $\mathfrak{f}_0=+\infty.$\\[0.2cm]
\hbox{}\quad
The following lemma gives us the precise value of $\mathfrak{f}_\xi.$ It easily follows from Lemma 2.3.\\[0.2cm]
\textbf{Lemma 2.5.} For $\xi\in F^{\times},$ we have
\[
\mathfrak{f}_\xi=\left\{\begin{array}{lr}
m-e+r \\ m-e
\end{array}
\begin{array}{lr}
\mbox{if }\xi=\varpi^{2m}u, u\in\mathcal{U}_r\setminus\mathcal{U}_{r+1},\\
\mbox{if }\ord(\xi)=2m+1.
\end{array}\right.
\]


\section{The Idempotents $e^K$ and $E^K$}
In this section we give two functions $e^K$ and $E^K$ on $SL_2(F)$ which are idempotents under certain convolution product.
They were introduced initially by Hiraga and Ikeda in \cite{IkedaHiraga} and play essential roles in the generalized Kohnen plus space.\\
\hbox{}\quad
Recall that $\Gamma=\Gamma[\mathfrak{d}^{-1},4\mathfrak{d}],$ we also set 
\[
\Lambda=\Gamma[\mathfrak{d}^{-1},\mathfrak{d}].
\]
\textbf{Definition 3.1.} The Hecke algebra $\widetilde{\mathcal{H}}=\widetilde{\mathcal{H}}(\widetilde{\Gamma}\widetilde{\backslash SL_2(F)}/\widetilde{\Gamma};\varepsilon)$ is defined to be the space of compactly supported locally constant function $\phi$ on
$\widetilde{SL_2(F)}$ such that $\phi(\gamma_1\widetilde{g}\gamma_2)=\varepsilon(\widetilde\gamma_1)\varepsilon(\widetilde\gamma_2)\phi(\widetilde{g})$ for any $\widetilde{g}\in SL_2(F),$ $\widetilde\gamma_1, \widetilde\gamma_2\in\Gamma.$ The multiplication of $\widetilde{\mathcal{H}}$ is given by the convolution product
\[
\phi_1\ast\phi_2=\int_{\widetilde{SL_2(F)}/\{\pm 1\}}\phi_1(\widetilde{g}\widetilde{h}^{-1})\phi_2(\widetilde{h})d\widetilde{h}.
\]
Here $d\widetilde{h}$ is taken to be the Haar measure on $\widetilde{SL_2(F)}$ such that the volume of $\widetilde{\Lambda}$ is $1$.\\[0.2cm]
\textbf{Definition 3.2.} The genius function $e^K$ on $\widetilde{SL_2(\mathfrak{o})}$ is defined as
\[
e^K(g)=\left\{\begin{array}{lr}
q^e(\phi_0,\omega_\psi(g)\phi_0)
\\
0
\end{array}
\begin{array}{lr}
\mbox{if }g\in\widetilde{\Lambda},
\\
\mbox{otherwise.}
\end{array}\right.
\]
\hbox{}\quad Recalling that $\varepsilon$ is given by $\omega_\psi(\widetilde{g})\phi_0=\varepsilon^{-1}(g)\phi_0$ on $\widetilde{\Gamma},$ we have $e^K\in\widetilde{\mathcal{H}}.$
The following lemma can be obtained by putting $\phi=\phi_0$ in equation (\ref{wc}) in Section 1. \\[0.2cm]
\textbf{Lemma 3.3.}\begin{itemize}
\item[(1)]For $a\in F^\times,$ we have
\[
\int_{t\in\mathfrak{o}}\psi(\frac{t^2}{a})dt=\alpha_\psi(a)\left|\frac{a}{2\bdelta}\right|^{1/2}\int_{t\in\mathfrak{o}}\psi(-\frac{at^2}{4\bdelta^2})dt
\]
\item[(2)]For nonzero $a\in \mathfrak{p}^{-c_\psi},$ we have
\[
\alpha_\psi(a)=|2a\bdelta|^{-1/2}\int_{t\in\mathfrak{o}}\psi(\frac{t^2}{4a\bdelta^2})dt
\]
\end{itemize}
\textbf{Lemma 3.4.} We have the following two properties about $e^K.$
\begin{itemize}
\item[(1)]
For nonzero $x\in\mathfrak{o},$ we have
\[
\begin{aligned}
e^K(\mathbf{u}^\flat(x\bdelta))
&=\alpha_\psi(x\bdelta)|2x|^{-1/2}\int_{t\in\mathfrak{o}}\overline{\psi(t^2/x\bdelta)}dt \\
&=|2|^{-1}\int_{t\in\mathfrak{o}}\psi(xt^2/4\bdelta)dt
\
\end{aligned}
\]
\item[(2)]
For $g=\left(\begin{array}{lr}a&b/\boldsymbol\delta \\ c\boldsymbol\delta&d\end{array}\right)\in\Lambda$ with $c\in\mathfrak{o}^\times,$ we have
\[
e^K([g])=\alpha_\psi(c\boldsymbol\delta)|2|^{-1/2}.
\]
\end{itemize}
\textit{Proof.}
Notice that
\[
\omega_\psi(\mathbf{u}^\flat(x\bdelta))\phi_0(s)=\overline{\alpha_\psi(x\bdelta)}|2x^{-1}|^{1/2}\int_{t\in F}\phi_0(t+s)\psi(t^2/x\bdelta)dt.
\]
So
\[
\begin{aligned}
(\phi_0,\omega_\psi(\mathbf{u}^\flat(x\bdelta))\phi_0)
&=\alpha_\psi(x\bdelta)|2x^{-1}|^{1/2}\int_{s\in F}\int_{t\in F}\phi_0(s)\phi_0(t+s)\overline{\psi(t^2/x\bdelta)}dtds \\
&=\alpha_\psi(x\bdelta)|2x^{-1}|^{1/2}\int_{t\in\mathfrak{o}}\overline{\psi(t^2/x\bdelta)}dt.
\end{aligned}
\]
Thus the first part of (1) follows. By using Lemma 3.3(1), we get the second part of (1). Now, (2) follows from the first part of (1).\done\\[0.2cm]
\textbf{Lemma 3.5.}
We have the following two identities:
\[
\begin{aligned}
&(1)\quad\int_{t\in\mathfrak{o}^\times}\alpha_\psi(t/\bdelta)dt=q^{-e/2}(1-q^{-1}).\\
&(2)\quad\int_{t\in\mathfrak{o}^\times}\alpha_\psi(\varpi t/\bdelta)dt=0.
\end{aligned}
\]
\textit{Proof.}
From Lemma 3.3(2), we have
\[
\int_{t\in\mathfrak{o}^\times}\alpha_\psi(t/\bdelta)dt=q^{e/2}\int_{x\in\mathfrak{o}}\int_{t\in\mathfrak{o}^\times}\psi\left(\frac{x^2t}{4\bdelta}\right)dtdx.
\]
Then (1) follows from
\begin{equation}
\label{intofkubota}
\int_{t\in\mathfrak{o}^\times}\psi(\xi t/\bdelta)dt=
\begin{cases}
1-q^{-1}&\mbox{ if }\xi\in\mathfrak{o},\\
-q^{-1}&\mbox{ if }\xi\in\varpi^{-1}\mathfrak{o}^\times,\\
0&\mbox{ otherwise}.
\end{cases}
\end{equation}
The proof of (2) is similar.\done\\[0.2cm]
\textbf{Lemma 3.6.}
If $\ord(x)$ is odd and $0<\ord(x)<2e$, then $e^K(\mathbf{u}^\flat(x\bdelta))=0.$\\[0.2cm]
\textit{Proof.} We may assume that the residue characteristic of $F$ is 2. Write $x=\varpi^{2r-1}u, u\in\mathfrak{o}^\times,$ where $0<r\leq e.$ Then
\[
\begin{aligned}
\int_{t\in\mathfrak{o}}\overline{\psi(t^2/x\bdelta)}dt&=\int_{t\in\mathfrak{o}}\int_{s\in\mathfrak{o}}\overline{\psi((t+\varpi^{r-1}s)^2/x\bdelta)}dsdt\\
&=\int_{t\in\mathfrak{o}}\int_{s\in\mathfrak{o}}\overline{\psi(t^2/x\bdelta)\psi(2\varpi^{r-1}st/x\bdelta)\psi(s^2/\varpi u\bdelta)}dsdt.
\end{aligned}
\]
Since $\ord(2\varpi^{r-1}/x)\geq 0,$ we have in the integral $\overline{\psi(2\varpi^{r-1}st/x\bdelta)}=1$. Moreover, because the mapping $s\rightarrow s^2$ is a bijection on $\mathfrak{o}/\mathfrak{p},$ we get that
\[
\int_{s\in\mathfrak{o}}\overline{\psi(s^2/\varpi u\bdelta)}ds=\int_{s\in\mathfrak{o}}\overline{\psi(s/\varpi u\bdelta)}ds=0.
\]
So then the lemma follows from (1) of Lemma 3.4.\done\\[0.2cm]
\textbf{Definition 3.7.}
The genuine function $E^K$ on $\widetilde{SL_2(F)}$ is defined as
\[
E^K(g)=e^K(\mathbf{w}_{2\bdelta}^{-1}g\mathbf{w}_{2\bdelta}).
\]
\quad
Noticing $\omega_\psi(\mathbf{w}_{2\bdelta}^{-1}\gamma\mathbf{w}_{2\bdelta})\phi_0=\omega_\psi(\gamma)\phi_0$ for $\gamma\in\widetilde{\Gamma},$ we get that $E^K\in\widetilde{\mathcal{H}}$ as $e^K$ does.
In fact, the conjugate mapping only changed the support of $e^K$ in its definition from $\widetilde{\Lambda}$ to $\widetilde{\Gamma[(4\mathfrak{d})^{-1},4\mathfrak{d}]}.$
That is 
\[
E^K(g)=\left\{\begin{array}{lr}
q^e(\phi_0,\omega_\psi(g)\phi_0)
\\
0
\end{array}
\begin{array}{lr}
\mbox{if }g\in\widetilde{\Gamma[(4\mathfrak{d})^{-1},4\mathfrak{d}]},
\\
\mbox{otherwise.}
\end{array}\right.
\]
\quad We denote the space of Schwartz functions on $\mathfrak{p}^{-e}/\mathfrak{o}$ by $\mathcal{S}(\mathfrak{p}^{-e}/\mathfrak{o}).$
That is, for any $\phi\in\mathcal{S}(\mathfrak{p}^{-e}/\mathfrak{o}),$ the support of $\phi$ is contained in $\mathfrak{p}^{-e}$ and $\phi(x+y)=\phi(x)$ for all $x\in\mathfrak{p}^{-e}, y\in\mathfrak{o}.$
It is easy to see that 
\[
\mathcal{S}(\mathfrak{p}^{-e}/\mathfrak{o})=\bigoplus_{\lambda\in\mathfrak{p}^{-e}/\mathfrak{o}}\mathbb{C}\cdot\phi_\lambda
\]
where $\phi_\lambda(x)=\phi_0(x-\lambda).$\\[0.2cm]
\textbf{Lemma 3.8.} The space $\mathcal{S}(\mathfrak{p}^{-e}/\mathfrak{o})$ forms an irreducible space under the action $\omega_\psi$ restricted to $\widetilde\Lambda.$\\[0.2cm]
\textit{Proof.}
In order to prove that $\mathcal{S}(\mathfrak{p}^{-e}/\mathfrak{o})$ is invariant under the action of $\widetilde{\Lambda}$, it suffices to show that for any $b\in\mathfrak{o}, c\in\mathfrak{o}^\times,$
$\omega_\psi(\mathbf{u}^\sharp(b/\bdelta))\phi, \omega_\psi(\mathbf{u}^\flat(b\bdelta))\phi, \omega_\psi(\mathbf{m}(c))\phi\in\mathcal{S}(\mathfrak{p}^{-e}/\mathfrak{o}),$
where $\phi$ is an arbitrary Schwartz function in $\mathcal{S}(\mathfrak{p}^{-e}/\mathfrak{o}).$
We omit the first and third cases since they are easy.
In the second case, it is easy to see $\omega_\psi(\mathbf{u}^\flat(b\bdelta))\phi(x+y)=\omega_\psi(\mathbf{u}^\flat(b\bdelta))\phi(x)$ where $x\in F, y\in\mathfrak{o}.$
Now if $x\notin\mathfrak{p}^{-e},$ we want to show that $\omega_\psi(\mathbf{u}^\flat(b\bdelta))\phi(x)=0.$
We may assume $\phi=\phi_\lambda$ for some $\lambda\in\mathfrak{p}^{-e}.$
Then
\[
\begin{aligned}
\omega_\psi(\mathbf{u}^\flat(b\bdelta))\phi(x)
&=\overline{\alpha_\psi(b\bdelta)}|2/b|^{1/2}\int_{t\in F}\phi_\lambda(x+t)\psi(t^2/b\bdelta)dt\\
&=\overline{\alpha_\psi(b\bdelta)}|2/b|^{1/2}\int_{t\in\mathfrak{o}}\psi((t-x+\lambda)^2/b\bdelta)dt.\\
\end{aligned}
\]
For simplicity, we denote $x-\lambda$ by $A,$ thus $\ord(A)=\ord(x)<-e.$
Setting $r=\ord(b)\geq0,$ the integral above becomes
\[
\begin{aligned}
\int_{t\in\mathfrak{o}}\psi((t-A)^2/b\bdelta)dt
&=\int_{s\in\mathfrak{o}}\int_{t\in\mathfrak{o}}\psi((t+\varpi^rs-A)^2/b\bdelta)dtds\\
&=\int_{t\in\mathfrak{o}}\psi((t-A)^2/b\bdelta)\int_{s\in\mathfrak{o}}\psi(2\varpi^rs(t-A)/b\bdelta)dsdt\\
&=0
\end{aligned}
\]
since $\ord(2\varpi^r(t-A)/b)<0.$
Hence $\omega_\psi(\mathbf{u}^\flat(b\bdelta))\phi(x)=0.$\\
\hbox{}\quad Now we prove that there is no non-trivial proper subspace of $\mathcal{S}(\mathfrak{p}^{-e}/\mathfrak{o})$ invariant under the action of $\widetilde{\Lambda}.$
 Let $\lambda, \mu\in\mathfrak{p}^{-e}$ such that $\lambda\not\equiv\mu$ (mod $\mathfrak{o}$).
We have
\[
\begin{aligned}
&\quad(\phi_\lambda,\omega_\psi(\mathbf{u}^\flat(z\bdelta))\phi_\mu)\\
&=\alpha_\psi(z\bdelta)|2z^{-1}|^{1/2}\int_{t\in F}\phi_0(t-\lambda)\int_{s\in F}\phi_0(t+s-\mu)\psi(-s^2/z\bdelta)dsdt\\
&=\alpha_\psi(z\bdelta)|2z^{-1}|^{1/2}\int_{t\in\mathfrak{o}}\int_{s\in\mathfrak{o}}\psi(-(s-t-\lambda+\mu)^2/z\bdelta)dsdt\\
&=\alpha_\psi(z\bdelta)|2z^{-1}|^{1/2}\int_{s\in\mathfrak{o}}\psi(-(s-\lambda+\mu)^2/z\bdelta)ds
\end{aligned}
\]
for any $z\in\mathfrak{o}.$\\
Now if we take $z$ to be in $\mathfrak{o}^\times,$ then the integral becomes
\[
\begin{aligned}
&\quad\int_{s\in\mathfrak{o}}\psi(-(s-\lambda+\mu)^2/z\bdelta)ds\\
&=\psi(-(\lambda-\mu)^2/z\bdelta)\int_{s\in\mathfrak{o}}\psi(-(s^2-2(\lambda-\mu))/z\bdelta)ds\\
&=\psi(-(\lambda-\mu)^2/z\bdelta)\\
&\neq0
\end{aligned}
\]
since $\ord(2(\lambda-\mu))\geq0.$\\
Thus we get that for some $\mathbf{u}^\flat(z\bdelta)\in\widetilde{\Lambda}, (\phi_\lambda,\omega_\psi(\mathbf{u}^\flat(z\bdelta))\phi_\mu)\neq0.$
The assertion we want to prove follows.\done
\\[0.2cm]
\textbf{Proposition 3.9.} The functions $e^K$ and $E^K$ are idempotents in $\widetilde{\mathcal{H}}.$\\[0.2cm]
\textit{Proof.}
Since the representation $\omega_\psi$ is unitary, the idempotence of $e^K$ follows from Lemma 3.8 and the Schur orthogonality relation.
Then according to the definition, we immediately get that $E^K$ is also idempotent.\done\\


\section{The Function $f^+_K$ and its Whittaker Function}
In this section we give the function $f^+_K$ on $\widetilde{SL_2(F)}.$ This function will later become the finite part of our Eisenstein series, and its Whittaker function is essential in the calculation of Fourier coefficients.
For more detail of both these functions, see \cite{IkedaHiraga}.\\ 
\hbox{}\quad
We set $B\subset SL_2(F)$ to be the Borel subgroup consisting of the upper triangular matrices and $\tilde{I}_\psi(s)$ to be the induced representation of $\widetilde{SL_2(F)}$ induced from the genuine character defined by
\[
\mathbf{u}^\sharp(b)\mathbf{m}(a)\longmapsto\frac{\alpha_\psi(1)}{\alpha_\psi(a)}|a|^s
\]
for $b\in F, a\in F^\times.$ Here $s$ is a complex number. Thus $\tilde{I}_\psi(s)$ is the space of genuine function $f$ on $\widetilde{SL_2(F)}$ satisfying
\[
f(\mathbf{u}^\sharp(b)\mathbf{m}(a)g)=\frac{\alpha_\psi(1)}{\alpha_\psi(a)}|a|^{s+1}f(g)
\]
for $b\in F, a\in F^\times.$
$\tilde{I}_\psi(s)$ can be considered as a genuine representation of $\widetilde{SL_2(F)}$ by the right translation. We let $q^{-2s}\neq q^{\pm 1},$ then it is well-known that under this restriction $\tilde{I}_\psi(s)$ is irreducible.\\[0.2cm]
\textbf{Definition 4.1.}
We set
\[
\mathcal{X}_r=
\left\{
\left(\begin{array}{lr}a&b/\bdelta\\c\bdelta&d\end{array}\right)\in\Lambda\,\bigg|\,
\ord(c)=2r
\right\}
\]
for $0\leq r<e$ and
\[
\mathcal{X}_e=\Gamma.
\]
\hbox{}\quad Notice that these sets are all bi-$\Gamma$-invariant subsets of $\Lambda.$ Also, we let $f_r\in\tilde{I}_\psi(s)$ be the function whose restriction to $\widetilde\Lambda$ is given by
\[
f_r(g)=\begin{cases}
q^{c_\psi s-e/2}\alpha_\psi(1)\overline{e^K(g)}&\mbox{if }g\in\widetilde{\mathcal{X}_r},\\
0 &\mbox{if }g\in\widetilde{\Lambda}\setminus\widetilde{\mathcal{X}_r}
\end{cases}
\] 
for $0\leq r\leq e.$\\[0.2cm]
\textbf{Definition 4.2.}
We define
\[
f_K^+=\sum^e_{r=0}q^{2sr-e+r}f_r,\qquad f^{[0]}_K=\sum^e_{r=0}f_r.
\]
\quad For any
$g=\left(\begin{array}{lr} a & b/\bdelta \\ c\bdelta & d \end{array}\right)\in\Lambda, b\neq 0$
since we have the decomposition
\[
[g]=\mathbf{u}^\sharp((a-1)/c\bdelta)\mathbf{u}^\flat(c\bdelta)\mathbf{u}^\sharp((d-1)/c\bdelta)
,\]
by Lemma 3.6, we get that the restriction of $f^{[0]}_K$ to $\widetilde\Lambda$ is
\[
f^{[0]}_K(g)=q^{c_\psi s-e/2}\alpha_\psi(1)\overline{e^K(g)}.
\]
So our definition of $f^{[0]}_K$ coincides with which in \textsection\,6 of \cite{IkedaHiraga}.\\[0.2cm]
\hbox{}\quad
The following proposition is an analogue of Proposition 4.4 in \cite{IkedaHiraga}, which shows that our $f_K^+$ is also the same as which in \textsection\,6 of \cite{IkedaHiraga}.\\[0.2cm]
\textbf{Proposition 4.3.} We have 
\[
\begin{array}{lr}
f_K^+(g\mathbf{w}_{2\bdelta})=\overline{\alpha_\psi(2\bdelta})q^{es-e/2}f^{[0]}_K(g),\\[0.2cm]
f_K^{[0]}(g\mathbf{w}_{2\bdelta})=\overline{\alpha_\psi(2\bdelta})q^{-es+e/2}f^+_K(g).
\end{array}
\]
\textit{Proof.}
It is easy to see that the two equations are equivalent, so we only have to prove the second one. It suffices to show that 
\[
f_{e-r}(g\mathbf{w}_{2\bdelta})=\overline{\alpha_\psi(2\bdelta)}q^{r-e/2-s(e-2r)}f_r(g)
\]
for $0\leq r\leq e.$ By the decomposition below Definition 4.2, we may assume that $g=\mathbf{u}^\flat(x\bdelta)$ where $x\in \mathfrak{o}$.
For $r=0,$ we have by Lemma 3.4(2),
\[
f_0(\mathbf{u}^\flat(x\bdelta))=\begin{cases}
q^{c_\psi s}\alpha_\psi(1)\overline{\alpha_\psi(x\bdelta)}&\mbox{if }\ord(x)=0,\\
0&\mbox{otherwise.}
\end{cases}
\]
On the other hand, 
\[
f_e(\mathbf{u}^\flat(x\bdelta)\mathbf{w}_{2\bdelta})=\begin{cases}
q^{c_\psi s-e/2-es}\alpha_\psi(1)\overline{\alpha_\psi(x\bdelta)}\overline{\alpha_\psi(2\bdelta)}&\mbox{if }\ord(x)=0,\\
0&\mbox{otherwise.}
\end{cases}
\]
Thus for the case $r=0$ the claim is true. The case $r=e$ is equivalent to the case $r=0.$ Now if $0<r<e,$ let $\ord(x)=2r$ then by the second part of Lemma 3.4(1), we have
\[
f_r(\mathbf{u}^\flat(x\bdelta))=q^{c_\psi s-e/2}\alpha_\psi(1)|2|^{-1}\int_{t\in\mathfrak{o}}\overline{\psi\left(\frac{xt^2}{4\bdelta}\right)}dt.
\]
On the other hand, since 
\[
\mathbf{u}^\flat(x\delta)\mathbf{w}_{2\bdelta}=[1_2,\langle-x\bdelta,-2x\rangle]\mathbf{u}^\sharp(1/x\bdelta)\mathbf{m}(-2/x)\mathbf{u}^\flat(-4\bdelta/x),
\]
by the first part of Lemma 3.4(1) we get that
\[
\begin{aligned}
f_{e-r}(\mathbf{u}^\flat(x\bdelta)\mathbf{w}_{2\bdelta})=&q^{c_\psi s-e/2}\alpha_\psi(1)\langle-x\bdelta,-2x\rangle\alpha_\psi(1)\alpha_\psi(2x)\left|\frac{2}{x}\right|^{s+1}\\
&\times\alpha_\psi(x\bdelta)\left|\frac{8}{x}\right|^{-1/2}\int_{t\in\mathfrak{o}}\overline{\psi\left(\frac{xt^2}{4\bdelta}\right)}dt.
\end{aligned}
\]
So through an easy calculation, the claim follows. The proof is completed.\done\\
\textbf{Definition 4.4.} For $\xi\in F^\times,\,g\in \widetilde{SL_2(F)},$ we define the Whittaker function $W^+_\xi(g)$ by
\[
W^+_{\xi}(g)=|\xi|^{1/2}\bgamma(\xi,q^{-s})^{-1}q^{\frak{f}_\xi s}\int_{x\in F}f^+_K(\mathbf{w}_1\mathbf{u}^\sharp(x)g)\overline{\psi(\xi x)}dx
\]
where
\[
\bgamma(\xi,y)=\frac{1-y^2/q}{1-\chi_\xi y/\sqrt{q}}.
\]
\quad The value of $W^+_\xi(1)$ is especially important for the calculation of the Fourier coefficients of the Eisenstein series in this paper.\\[0.2cm]
\textbf{Proposition 4.5.} (\textsection\,6 of \cite{IkedaHiraga}) We have
\[
W^+_\xi(\mathbf{m}(a))=\frac{\alpha_\psi(1)}{\alpha_\psi(a)}|\xi a^2|^{1/2}\Psi(\xi a^2,q^{-s}),\quad a\in F^\times. 
\]
In particular,
\[
W^+_\xi(1)=|\xi|^{1/2}\Psi(\xi,q^{-s}).
\]
Here
\[
\large{
\Psi(\xi,y)=\begin{cases}
\frac{y^{\mathfrak{f}_\xi+1}-y^{-\mathfrak{f}_\xi-1}}{y-y^{-1}}-\chi_\xi q^{-1/2}\frac{y^{\mathfrak{f}_\xi}-x^{-\mathfrak{f}_\xi}}{y-y^{-1}} &\mbox{if }\mathfrak{f}_\xi\geq 0, \\
0 &\mbox{if }\mathfrak{f}_\xi< 0.
\end{cases}
}
\]
\hbox{}\quad The Whittaker function is not defined for the case $\xi=0.$
However, we can get the explicit value of a similar integral for that case.\\[0.2cm]
\textbf{Lemma 4.6.} 
Assume that $e>0.$
Then we have
\[
\int_{t\in F}f_r(\mathbf{w}_1\mathbf{u}^\sharp(t))dt=
\begin{cases}
1&\mbox{ if }r=0,\\
(1-q^{-1})q^{-2rs}&\mbox{ if }0<r<e,\\
(1-q^{-1}){q^{-2es}}/(1-q^{-2s})&\mbox{ if }r=e.
\end{cases}
\]
\textit{Proof.}
The case $r=0$ follows from that $f_0(\mathbf{w}_1\mathbf{u}^\sharp(t))=q^{-c_\psi}$ if $\ord(t)\geq-c_\psi$ by Lemma 3.4(2) and vanishes otherwise.
Next, we consider the case $r=e.$
Notice that $f_e(\mathbf{w}_1\mathbf{u}^\sharp(t))=q^{c_\psi s+e/2}\alpha_\psi(t)|t|^{-s-1}$ if $\ord(t)\leq-2e-c_\psi$ by Lemma 1.1 and vanishes otherwise.
So
\[
\begin{aligned}
\int_{t\in F}f_e(\mathbf{w}_1\mathbf{u}^\sharp(t))dt=
&q^{c_\psi s+e/2}\sum_{n=2e+c_\psi}^\infty\int_{t\in\varpi^{-n}\mathfrak{o}^\times}\alpha_\psi(t)|t|^{-s-1}dt\\
=&q^{e/2}\sum_{n=2e}^\infty q^{-ns}\int_{t\in\mathfrak{o}^\times}\alpha_\psi(\varpi^{-n}t/\bdelta)dt\\
=&q^{e/2-2se}(1-q^{-2s})^{-1}\left(\int_{t\in\mathfrak{o}^\times}\alpha(t/\bdelta)dt+q^{-s}\int_{t\in\mathfrak{o}^\times}\alpha(\varpi t/\bdelta)dt\right).
\end{aligned}
\]
Then the assertion for $r=e$ follows from Lemma 3.5.
Finally, let $0<r<e$.
Notice that by the first part of Lemma 3.4(1), we have 
\[
f_r(\mathbf{w}_1\mathbf{u}^\sharp(t))=
\begin{cases}
q^{c_\psi(s-1/2)}|t|^{-s-1/2}\int_{x\in\mathfrak{o}}\psi(tx^2)dx&\mbox{ if }\ord(t)=-c_\psi-2r,\\
0&\mbox{ otherwise.}
\end{cases}
\]
It follows that
\[
\begin{aligned}
&\int_{t\in F}f_r(\mathbf{w}_1\mathbf{u}^\sharp(t))dt\\
=&q^{c_\psi(s-1/2)}\int_{t\in\varpi^{-c_\psi-2r}\mathfrak{o}^\times}|t|^{-s-1/2}\int_{x\in\mathfrak{o}}\psi(tx^2)dxdt\\
=&q^{r(1-2s)}\int_{x\in\mathfrak{o}}\int_{t\in\mathfrak{o}^\times}\psi\left(\frac{tx^2}{\varpi^{2r}\bdelta}\right)dtdx.
\end{aligned}
\]
Now using the identity (\ref{intofkubota}) in the proof of Lemma 3.5, we get the assertion for $0<r<e.$\done\\[0.2cm]
\textbf{Lemma 4.7.} Assume that $e=0.$ Then we have
\[
\int_{t\in F}f_e(\mathbf{w}_1\mathbf{u}^\sharp(t))dt=\frac{1-q^{-2s-1}}{1-q^{-2s}}.
\]
\textit{Proof.} Similar to the proof of the case $r=e$ in the last lemma.\done\\[0.2cm]
From the two lemmas and the definition of $f^+_K,$ we obtain the following proposition.\\[0.2cm]
\textbf{Proposition 4.8.} For $a\in F^\times,$ we have
\[
\int_{t\in F}f^+_K(\mathbf{w}_1\mathbf{u}^\sharp(t)\mathbf{m}(a))dt=\frac{\alpha_\psi(1)}{\alpha_\psi(a)}|a|^{-s+1}\frac{1-q^{-2s-1}}{1-q^{-2s}}.
\]
In particular,
\[
\int_{t\in F}f^+_K(\mathbf{w}_1\mathbf{u}^\sharp(t))dt=\frac{1-q^{-2s-1}}{1-q^{-2s}}.
\]
\quad Recall that $\tilde{I}_\psi(s)$ can be considered as a genuine representation of $\widetilde{SL_2(F)}$ by right translation $\rho$, so the Hecke algebra $\widetilde{\mathcal{H}}$ acts on $\tilde{i}_\psi(s)$ by
\[
\rho(\phi)f(g)=\int_{\widetilde{SL_2(F)}/\{\pm 1\}}\phi(h)f(gh)dh
\]
for $\phi\in\widetilde{\mathcal{H}}, f\in\tilde{i}_\psi(s)$.
In particular, $\rho(e^K)f^{[0]}_K=f^{[0]}_K\ast\overline{e^K}=f^{[0]}_K$ by the idempotence of $e^K.$
Similarly, $\rho(E^K)f^+_K=f^+_K.$\\
\hbox{}\quad The following proposition reveals some nature of the structure of the Kohnen plus space.\\[0.2cm]
\textbf{Theorem 4.9.}
If we set $\tilde{I}_\psi(s)^{e^K}=\{f\,|\,f\in\tilde{I}_\psi(s),\,\rho(e^K)f=f\}$ and $\tilde{I}_\psi(s)^{E^K}=\{f\,|\,f\in\tilde{I}_\psi(s),\,\rho(E^K)f=f\},$
then $\tilde{I}_\psi(s)^{e^K}=\mathbb{C}\cdot f^{[0]}_K$ and $\tilde{I}_\psi(s)^{E^K}=\mathbb{C}\cdot f^+_K.$\\[0.2cm]
\textit{Proof.}
The first identity follows from that any $f\in\tilde{I}_\psi(s)$ is determined by its restriction to $\widetilde{\Lambda}.$
Then the second identity follows from the statement that $f\in\tilde{I}_\psi(s)^{E^K}$ if and only if $\rho(\mathbf{w}_{2\bdelta})f\in\tilde{I}_\psi(s)^{e^K}$ which is easy to check.\done


\section{The Archimedean Case}
This section refers to \textsection\,8 of \cite{IkedaHiraga}. In this section, we set $F=\mathbb{R}$ and $\psi(x)=\mathbf{e}(x)$ for $x\in F.$
The Weil constant is $\alpha_\psi(a)$ is given by
\[
\alpha_\psi(a)=
\begin{cases}
\exp(\pi\sqrt{-1}/4)&\mbox{if }a>0,\\
\exp(-\pi\sqrt{-1}/4)&\mbox{if }a<0.
\end{cases}
\]
\quad The real metaplectic group $\widetilde{SL_2(\mathbb{R})}$ is the unique non-trivial topological double covering of $SL_2(\mathbb{R})$ equipped with the multiplication mentioned in Section 1.
It acts on the upper-half plane $\mathfrak{h}$ through $SL_2(\mathbb{R})$ by the usual action.\\
\hbox{}\quad We define a factor of automorphy $\tilde{j} : \widetilde{SL_2(\mathbb{R})}\times\mathfrak{h}\rightarrow\mathbb{C}$ by
\[
\tilde{j}\left(\left[
\begin{pmatrix}
a&b\\c&d
\end{pmatrix}
,\zeta\right],\tau\right)=
\begin{cases}
\zeta\sqrt{d}&\mbox{if }c=0, d>0,\\
-\zeta\sqrt{d}&\mbox{if }c=0, d<0,\\
\zeta(c\tau+d)^{1/2}&\mbox{if }c\neq 0.
\end{cases}
\]
Then $\tilde{j}$ is the unique factor of automorphy such that $\tilde{j}^2=j,$ where $j$ is the usual factor of automorphy on $SL_2(\mathbb{R})\times\mathfrak{h}.$


\section{Automorphic Forms on $\widetilde{SL_2(\mathbb{A})}$}
In this section we give a quick introduction of automorphic forms on $\widetilde{SL_2(\mathbb{A})}.$
Only those definitions and properties we need will be stated here.
For more detail, one can consult \textsection\,8 of \cite{IkedaHiraga}.\\
\hbox{}\quad We set $F$ to be a totally real number field with degree $n$ over $Q$ and denote its adele ring, ring of integers and different over $Q$ by $\mathbb{A},$ $\mathfrak{o}$ and $\mathfrak{d}.$
The $n$ infinite places of $F$ are denoted by $\infty_1,\infty_2,...,\infty_n.$
Let $\psi=\prod_{v}\psi_v$ be a non-trivial addictive character of $\mathbb{A}/F,$ here $v$ runs over all finite and infinite places of $F$.
If $v$ is a real place, $\psi_v$ is given by $\psi_v(x)=\mathbf{e}(x).$
For any finite place $v$, we assume that, as in Section 1, the order of  $\psi_v$ is $c_v,$ where $c_v$ is the multiplication of the corresponding prime ideal $\mathfrak{p}_v$ in the prime decomposition of $\mathfrak{d}$.
From now, the lower subscript $v$ simply means the relative element is with respect to the local field $F_v.$\\
\hbox{}\quad Let us construct the metaplectic double covering of $SL_2(\mathbb{A}).$
We use the method Hiraga and Ikeda used in \cite{IkedaHiraga}.
Let $\mathfrak{S}$ be a set of bad places of $F,$ which contains all places $v$ which is even or infinite, or such that $c_v\neq0.$
Set
\[
SL_2(\mathbb{A})_\mathfrak{S}=\prod_{v\in\mathfrak{S}}SL_2(F_v)\times\prod_{v\notin\mathfrak{S}}SL_2(\mathfrak{o}_v).
\] 
The metaplectic double covering $\widetilde{SL_2(\mathbb{A})_\mathfrak{S}}$ of $SL_2(\mathbb{A})_\mathfrak{S}$ is defined by the 2-cocycle $\prod_{v\in\mathfrak{S}}\mathbf{c}_v(g_{1,v},g_{2,v}).$\\
\hbox{}\quad For $\mathfrak{S}\subset\mathfrak{S}',$ we define an embedding
\[
\widetilde{SL_2(\mathbb{A})_\mathfrak{S}}\rightarrow\widetilde{SL_2(\mathbb{A})_\mathfrak{S'}}
\]
by
\[
[(g_v),\zeta]\mapsto[(g_v),\zeta\prod_{v\in\mathfrak{S}'\setminus\mathfrak{S}\mathbf{s}_v(g_v)}\mathbf{s}_v(g_v)]
\]
where $\mathbf{s}_v:SL_2(\mathfrak{o}_v)\rightarrow\widetilde{SL_2(\mathfrak{o}_v)}$ is the unique splitting of the covering $\widetilde{SL_2(\mathfrak{o}_v)}\rightarrow SL_2(\mathfrak{o}_v)$ for $v\notin\mathfrak{S}.$
Now the metaplectic double covering of $\widetilde{SL_2(\mathbb{A})}$ is defined to be the direct limit $lim_{\rightarrow}\widetilde{SL_2(\mathbb{A})_\mathfrak{S}}$
and there exists a canonical embedding $\widetilde{SL_2(F_v)}\hookrightarrow\widetilde{SL_2(\mathbb{A})}$ for each place $v.$
It is well-known that $SL_2(F)$ can be canonically embedded into $\widetilde{SL_2(\mathbb{A})}$ and that the embedding is given by $\gamma\mapsto[\gamma,1]$ for sufficiently large $\mathfrak{S}.$
From now we see $SL_2(F)$ as a subgroup of $\widetilde{SL_2(\mathbb{A})}$ through this embedding.
If $\prod'_v\widetilde{SL_2(F_v)}$ is the restricted product with respect to $\mathbf{s}_v(SL_2(\mathfrak{o}_v))$ then there is a canonical subjection $\prod'_v\widetilde{SL_2(F_v)}\rightarrow\widetilde{SL_2(\mathbb{A})}.$
We denote the image of $(g_v)$ also by $(g_v).$
As in the local case, for $x\in\mathbb{A}$ and $a\in\mathbb{A}^\times,$ we let
\[
\begin{aligned}
&\mathbf{u}^\sharp(x)=(\mathbf{u}^\sharp(x_v)),\quad\mathbf{u}^\flat(x)=(\mathbf{u}^\flat(x_v)),\\
&\mathbf{m}(a)=(\mathbf{m}(a_v)),\quad\mathbf{w}_a=(\mathbf{w}_{a_v}).
\end{aligned}
\]
\hbox{}\quad A function $f$ on $\widetilde{SL_2(\mathbb{A})}$ is said to be genuine if $f([g,\zeta])=\zeta f([g,1])$ for $g\in SL_2(\mathbb{A}), \zeta\in\{\pm1\}.$
Given a family of genuine functions $f_v$ for each place $v$ such that $f_v(g_v)=1$ for almost all $v$ where $g_v\in\mathbf{s}_v(\mathfrak{o}_v),$
the product $\prod_vf_v$ defined by $(\prod_vf_v)(g)=\prod_vf_v(g_v)$ gives a genuine function on  $\widetilde{SL_2(\mathbb{A})}.$\\
\hbox{}\quad Recall that $\Gamma=\Gamma[\mathfrak{d}^{-1},4\mathfrak{d}]$ and $\Gamma_v=\Gamma[\mathfrak{d}_v^{-1},4\mathfrak{d}_v].$
The restricted product $\Gamma'_f=\prod'_{v<\infty}\Gamma_v$ with respected to $SL_2(\mathfrak{o_v})$ for each $v$ is a compact open subgroup of $SL_2(\mathbb{A}_f),$ where $\mathbb{A}_f$ is the finite part of $\mathbb{A}.$
From the last paragraph of Section 1 there exists a genuine character $\varepsilon_v$ on $\widetilde{\Gamma_v}$ such that $\omega_{\psi_v}(g_v)\phi_{0,v}=\varepsilon_v(g_v)^{-1}\phi_{0,v}$ for each finite place $v.$
We denote $\varepsilon=\prod_{v<\infty}\varepsilon_v$ the genuine character on $\widetilde{\Gamma_f'}.$
Let $\kappa=(\kappa_1,\kappa_2,...,\kappa_n)\in\mathbb{Z}^n$ be an $n$-tuple with non-negative integers as whose components.
We define a factor of automorphy $j^{\kappa+1/2}(\gamma,z)$ for $\gamma\in\Gamma$ and $z=(z_1,z_2,...,z_n)\in\mathfrak{h}^n$ by
\[
j^{\kappa+1/2}(\gamma,z)=\prod_{v<\infty}\varepsilon_v([\gamma,1])\prod_{i=1}^{n}\tilde{j}([\iota_i(\gamma),1],z_i)^{2\kappa_i+1}.
\]
\textbf{Definition 6.1.} We denote $M_{\kappa+1/2}(\Gamma)$ and $S_{\kappa+1/2}(\Gamma)$ the spaces of Hilbert modular forms and cusp forms for $\Gamma$ of weight $\kappa+1/2$ with respect to the factor of automorphy $j^{k+1/2}(\gamma,z).$\\[0.2cm]
\hbox{}\quad For any $h\in M_{\kappa+1/2}(\Gamma),$ we can consider it as an automorphic form on $SL_2(F)\backslash\widetilde{SL_2(\mathbb{A})}$ as follows.
For any $\tilde{g}\in\widetilde{SL_2(\mathbb{A})},$ by the strong approximation theorem for $SL_2(\mathbb{A})$ there exist $\gamma\in SL_2(F),\,\tilde{g}_\infty\in\widetilde{SL_2(\mathbb{R})^n}$ and $\tilde{g}_f\in\widetilde{\Gamma_f}$ such that $\tilde{g}=\gamma\tilde{g}_\infty\tilde{g}_f.$
Then if we put
\begin{equation}
\label{modutoauto}
\varphi_h(\tilde{g})=h(\tilde{g}_\infty(\mathbf{i}))\varepsilon(\tilde{g}_f)^{-1}\prod_{i=1}^n(\tilde{j}(\tilde{g}_{\infty_i},\mathbf{i})^{2\kappa_i+1})^{-1}
\end{equation}
where $\mathbf{i}=(\sqrt{-1},...,\sqrt{-1})\in\mathfrak{h}^n,$ $\varphi_h$ is a genuine automorphic form on $SL_2(F)\backslash\widetilde{SL_2(\mathbb{A})}$.
We set
\[
\begin{array}{lr}
\mathbf{A}_{\kappa+1/2}(SL_2(F)\backslash\widetilde{SL_2(\mathbb{A})}; \widetilde{\Gamma_f'},\varepsilon)=\{\varphi_h\,|\,h\in M_{\kappa+1/2}(\Gamma)\},\\[0.2cm]
\mathbf{A}_{\kappa+1/2}^{\mbox{cusp}}(SL_2(F)\backslash\widetilde{SL_2(\mathbb{A})}; \widetilde{\Gamma_f'},\varepsilon)=\{\varphi_h\,|\,h\in S_{\kappa+1/2}(\Gamma)\}.
\end{array}
\]
\quad On the contrary, for each $\varphi\in\mathbf{A}_{\kappa+1/2}(SL_2(F)\backslash\widetilde{SL_2(\mathbb{A})}; \widetilde{\Gamma'_f}, \varepsilon),$ we can refer it to an Hilbert modular form $h\in M_{\kappa+1/2}(\Gamma)$ as follows.
For any $z\in\mathfrak{h}^n,$ there exist $\tilde{g}_\infty=[(g_{\infty_1},...,g_{\infty_n}),1]\in\widetilde{SL_2(\mathbb{R})}$ such that $z=\tilde{g}_\infty(\mathbf{i}).$
Then we define $h$ by
\begin{equation}
\label{autotomodu}
h(z)=\varphi(\tilde{g}_\infty)\prod_{i=1}^n\tilde{j}(\tilde{g}_{\infty_i},\mathbf{i})^{2\kappa_i+1}.
\end{equation}
It is obvious that $\varphi=\varphi_h.$
Moreover, for any adele $h\in\mathbb{A},\,\tilde{g}\in\widetilde{SL_2(\mathbb{A})},$ we have
\[
\varphi(\mathbf{u}^\sharp(h+x)\tilde{g})=\varphi(\mathbf{u}^\sharp(h)\tilde{g})
\]
where $x\in F$ is arbitrarily chosen.
So taking $dh$ as the Haar measure of $\mathbb{A}/F$ such that Vol$(\mathbb{A}/F)=1,$ we get
\begin{equation}
\label{expansion22}
\varphi(\mathbf{u}^\sharp(x)\tilde{g})=\sum_{\xi\in F}c(\xi,\tilde{g})\psi(\xi x)
\end{equation}
where
\[
c(\xi,\tilde{g})=\int_{\mathbb{A}/F}\varphi(\mathbf{u}^\sharp(h)\tilde{g})\overline{\psi(\xi x)}dh.
\]
On the other hand, if $h$ has the Fourier expansion $h(z)=\sum_{\xi\in F}c(\xi)q^\xi,$ from equation (\ref{modutoauto}) we find
\[
\varphi(\tilde{g})=\sum_{\xi\in F}c(\xi)q^\xi\varepsilon(\tilde{g}_f)^{-1}\prod_{i=1}^n\tilde{j}(\tilde{g}_{\infty_i},\mathbf{i})^{-2\kappa_i-1}.
\]
If now $g_f=1$ and $x_\infty\in\mathbb{A}$ is some adele such that $x_{\infty,v}=0$ for all $v<\infty,$ we have
\[
\varphi(\mathbf{u}^\sharp(x_\infty)\tilde{g})=\sum_{\xi\in F}c(\xi)q^\xi\textbf{e}(\sum_i^n\iota_i(\xi)x_{\infty_i})\prod_{i=1}^n\tilde{j}(\tilde{g}_{\infty_i},\mathbf{i})^{-2\kappa_i-1}.
\]
Comparing this with equation (\ref{expansion22}) above, we get the relation between $c(\xi,\tilde{g})$ and $c(\xi).$\\[0.2cm]
\textbf{Lemma 6.2.} let $z=\tilde{g}_\infty(\mathbf{i}).$ The following identity holds:
\[
c(\xi,\tilde{g_\infty})=c(\xi)q^\xi\prod_{i=1}^n\tilde{j}(\tilde{g}_{\infty_i},\mathbf{i})^{-2\kappa_i-1}.
\]
\hbox{}\quad Recall that $\Lambda_v=\Gamma[\mathfrak{d}_v^{-1},\mathfrak{d}_v]$ for $v<\infty.$
The restricted product $\prod'_{v<\infty}\Lambda_v$ is a maximal compact subgroup of $SL_2(\mathbb{A}_f).$
We denote it by $\Lambda_f'.$
If $\Omega'_{\psi,v}$ is the restriction of the Weil representation $\omega_{\psi,v}$ to $\widetilde{\Lambda_v},$
then the subrepresentation $\Omega_{\psi,v}$ of $\Omega'_{\psi,v}$ on $\mathcal{S}(2^{-1}\mathfrak{o}_v/\mathfrak{o}_v)$ is irreducible by Lemma 3.8.
So we get an irreducible representation $\Omega_\psi=\prod_{v<\infty}\Omega_{\psi,v}$ of $\widetilde{\Lambda_f'}$ on $\mathcal{S}(2^{-1}\hat{\mathfrak{o}}/\hat{\mathfrak{o}}),$ where $\hat{\mathfrak{o}}=\prod_{v<\infty}\mathfrak{o}_v.$
Note that $2^{-1}\hat{\mathfrak{o}}/\hat{\mathfrak{o}}$ can be canonically identified with $2^{-1}\mathfrak{o}/\mathfrak{o}.$


\section{The Kohnen Plus Space}
The concept of Kohnen plus space was initially introduced by Kohnen in 1980 \cite{Kohnen:80}.
It is defined as the subspace of $M_{\kappa+1/2}(4)$ consisting of modular forms whose $n$-th Fourier coefficients vanish whenever $(-1)^\kappa n\equiv 2,3$ (mod 4),
where $\kappa$ is an integer $\geq 3$ and $M_{\kappa+1/2}(4)$ denotes the space of modular forms of weight $\kappa+1/2$ on the congruence group $\Gamma_0(4).$
We denote the subspace by $M_{\kappa+1/2}^+(4).$
We also denote $M^+_{\kappa+1/2}(4)\cap S_{\kappa+1/2}(4)$ by $S^+_{\kappa+1/2}(4).$
Here $S_{\kappa+1/2}(4)$ is the subspace of cusp forms lying in $M_{\kappa+1/2}(4).$
The operators $U_4$ and $W_4$ on $M_{\kappa+1/2}(4)$ are defined by
\[
({U_4f})(z)=\frac{1}{4}\sum_{l\mbox{ mod }4}f\left(\frac{z+l}{4}\right),
\]\\[-0.8cm]
\[
(W_4f)(z)=(-2iz)^{-\kappa-1/2}f\left(-\frac{1}{4z}\right).
\]
Kohnen \cite{Kohnen:80} proved that $S_{\kappa+1/2}^+(4)$ is the eigenspace of $W_4U_4$ with respect to the eigenvalue $(-1)^{\kappa(\kappa+1)/2}2^\kappa.$\\
\hbox{}\quad Recently, Hiraga and Ikeda \cite{IkedaHiraga} gave a generalization of the plus space to the case for usual Hilbert modular forms.
That is what we want to introduce in this section.\\
\hbox{}\quad We use the same notations as in the last section.
Also, we fix an unit $\eta\in\mathfrak{o}^\times$ such that $N_{F/\mathbb{Q}}(\eta)=\prod_{i=1}^n(-1)^{\kappa_i}.$
We only consider the condition that such an unit exists.\\[0.2cm]
\textbf{Definition 7.1.}
For each $\xi\in F,$ we write $\xi\equiv\square$ mod 4 if there exists $y\in\mathfrak{o}$ such that $\xi\equiv y^2$ mod $4\mathfrak{o}.$\\[0.2cm]
\hbox{}\quad For any $\xi\in F$ and a non-archimedean place $v$ of $F,$ let $\mathfrak{f}_{\xi,v}$ be the invariant as which defined in Definition 2.4 with respect to the local field $F_v.$
Then we have the following lemma.\\[0.2cm]
\textbf{Lemma 7.2.}
(Lemma 13.1 of \cite{IkedaHiraga})
For each $\xi\in F,$ $\xi\equiv\square$ mod 4 if and only if $\mathfrak{f}_{\xi,v}\geq 0$ for any finite place $v$ of $F$.\\[0.2cm]
\hbox{}\quad Notice that if $\xi\equiv\square$ mod 4, then $\xi\in\mathfrak{o}.$\\[0.2cm]
\textbf{Definition 7.3.}
The Kohnen plus space $M^+_{\kappa+1/2}(\Gamma)$ is the space of $h(z)=\sum_{\xi\in\mathfrak{o}}c(\xi)\mathbf{e}(\xi z)\in M_{\kappa+1/2}(\Gamma)$ such that $c(\xi)$ vanishes unless $\eta\xi\equiv\square$ mod 4.
We also put $S^+_{\kappa+1/2}(\Gamma)=M^+_{\kappa+1/2}(\Gamma)\cap S_{\kappa+1/2}(\Gamma),$ which is also called a Kohnen plus space.\\[0.2cm]
\hbox{}\quad For any finite place $v,$ let $E^K_v$ be the idempotent defined in Section 3 with respect to $F_v,$
then $\prod_{v<\infty}E^K_v,$ which we denote by $E^K$ in the global case, is a genuine locally constant function on $\widetilde{SL_2}(\mathbb{A}_2)$ with compact support.
$E^K$ naturally acts on $M_{\kappa+1/2}(\Gamma)$ at each finite place by $\rho$ as in Section 4.
Let $M_{\kappa+1/2}(\Gamma)^{E^K}=\{h\in M_{\kappa+1/2}(\Gamma)\,|\,\rho(E^K)h=h\}$ and $S_{\kappa+1/2}(\Gamma)^{E^K}=M^{E^K}_{\kappa+1/2}(\Gamma)\cap S_{\kappa+1/2}(\Gamma).$
Hiraga and Ikeda showed that the two spaces coincide with the two plus spaces we just defined.\\[0.2cm]
\textbf{Theorem 7.4.}
(Theorem 13.5 of \cite{IkedaHiraga})
We have
\[
M_{\kappa+1/2}(\Gamma)^{E^K}=M^+_{\kappa+1/2}(\Gamma),\qquad S_{\kappa+1/2}(\Gamma)^{E^K}=S^+_{\kappa+1/2}(\Gamma).
\]


\section{Construction of the Eisenstein Series}
Now we can begin to construct our Eisenstein series.
Again, we use the same notations as in the last section, but now the $n$-tuple $\kappa$ is restricted to be parallel and positive, that is, $\kappa=(k,...,k)$ for certain $k\geq 1.$
We simply denote the component $k$ of $\kappa$ also by $\kappa$ if there is no confusion.
The unit $\eta$ is chosen to be $(-1)^\kappa\in\mathfrak{o}^\times.$ 
Furthermore, for simplicity, we let $G=SL_2(F)$ and $B$ be the Borel subgroup of $G$ consisting of all upper triangular matrices.
By Theorem 7.4 and Theorem 4.9, the finite part of the Eisenstein series must has the form
\[
\prod_{v<\infty}f^+_{K,v}
\]
where $f^+_{K,v}\in\tilde{I}_{\psi,v}(s_v)\mbox{ for some }s_v\in\mathbb{C}.$
But instead of $\psi$ we considered so far, let us use the new non-trivial additive character $\psi_\eta$ defined by
\[
\psi_\eta(a)=\psi(\eta a).
\]
Obviously $\psi_{\eta,v}$ has the same order with $\psi_v$ for each finite place $v.$
The only things we have to be careful are that now the Weil constant of any $a\in F_v^\times$ with respect to $\psi_{\eta,v}$ is $\alpha_{\psi_v}(\eta a)$ for each finite and infinite place $v$
and with respect to each infinite place the character is now $\mathbf{e}(\eta x).$
Now suppose that $F$ has (wide) class number $h$ and ideal class group $\mathcal{C_F}.$
There exist $h$ class characters $\chi_1=1,\chi_2,...,\chi_h$ of $\mathcal{C}_F$.
Each character $\chi_j$ induces a Hecke character of $F,$ which we also denote by $\chi_j.$
We fix one character $\chi'=\chi_j$ for some $1\leq j\leq h.$
Then
\[
\chi'_{\infty_i}=1
\]
for $1\leq i\leq n$ and
\[
\chi'_v(a)=|a|_v^{\mathfrak{s}_v}, a\in F_v^\times
\]
where $\mathfrak{s}_v$ is some complex number, for each $v<\infty.$\\
Now let
\[
s_v=\kappa-1/2+\mathfrak{s}_v
\]
for each $v<\infty.$\\
\hbox{}\quad Assume the condition $\kappa\geq2$.
We define a function $f=f_{\chi'}=\prod_{v}f_v$ on $\widetilde{SL_2(\mathbb{A})}$ by
\begin{equation}
\label{kgeq2}
f_v(g)=\begin{cases}
f^+_{K,v}(g)\in\tilde{I}_{\psi_\eta,v}(s_v)&\mbox{ if } v<\infty,\\
\tilde{j}(g_v,\mathbf{i})^{-(2\kappa+1)}&\mbox{ if } v\,|\,\infty.
\end{cases}
\end{equation}
\hbox{}\quad One can check that $f$ is left-$B$-invariant using the well-known fact that
$
\prod_{v\leq\infty}\alpha_{\psi_v}(a)=1
$
for any $a\in F^\times$(and this is the reason why we use $\psi_\eta$ instead of $\psi$).
Now our Eisenstein series is defined as follows.\\[0.2cm]
\textbf{Definition 8.1.}
For $\kappa\geq2,$ we define the Eisenstein series $E'(g)=E'_{\kappa+1/2}(g,\chi')$ on $\widetilde{SL_2(\mathbb{A})}$ of weight $\kappa+1/2$ twisted by $\chi'$ as
\[
E'(g)=\sum_{\gamma\in B\backslash G}f(\gamma g)
\]
where $g\in\widetilde{SL_2(\mathbb{A})}.$\\[0.2cm]
\hbox{}\quad Since the series in the definition does not converge if we set $\kappa=1.$
We should use another way to define the Eisenstein series of weight 3/2.
Now let $F\neq\mathbb{Q}$ and $\kappa=1.$
So for each finite place $v,$ $s_v=1/2+\mathfrak{s}_v.$
Instead of the function $f$ defined by equation (\ref{kgeq2}), for any $\epsilon\in\mathbb{C}$ such that Re$(\epsilon)>1/2,$ we set a function $f_\epsilon=f_{\chi',\epsilon}=\prod_v f_{\epsilon,v}$ on $\widetilde{SL_2(\mathbb{A})}$ by
\begin{equation}
\label{keq1}
f_{\epsilon,v}(g)=\begin{cases}
f^+_{K,v}(g)\in\tilde{I}_{\psi_\eta,v}(s_v+\epsilon)&\mbox{ if } v<\infty,\\
\tilde{j}(g_v,\mathbf{i})^{-3}|\tilde{j}(g_v,\mathbf{i})|^{-2\epsilon}&\mbox{ if } v\,|\,\infty.
\end{cases}
\end{equation}
Again one notice that $f_\epsilon$ is left-$B$-invariant.
For $g\in\widetilde{SL_2(\mathbb{A})},$ the function $E_\epsilon(g)'$ is defined  by
\[
E'_\epsilon(g)=\sum_{\gamma\in B\backslash G}f_\epsilon(\gamma g)
\]
converges and can be continued analytically into the whole $\epsilon-$plane.
So the definition of Eisenstein series of weight 3/2 is as following.\\[0.2cm]
\textbf{Definition 8.2.}
If $F\neq\mathbb{Q},$ we define the Eisenstein series $E'(g)=E'_{3/2}(g,\chi')$ on $\widetilde{SL_2(\mathbb{A})}$ of weight $3/2$ twisted by $\chi'$ as
\[
E'(g)=E'_0(g)
\]
where $g\in\widetilde{SL_2(\mathbb{A})}$ and $E'_\epsilon$ is the series defined as above.\\[0.2cm]
\textit{Remark.}
It is known that the modular form $\mathcal{H}_\kappa$ constructed by Cohen, which our Eisenstein series will coincide with for $F=\mathbb{Q}$ and $\kappa\geq2,$ is not a modular form if $\kappa=1.$
We will also see that if $F=\mathbb{Q},$ the series $E'_{3/2}(g,\chi')$ defined as above does not give a modular form.
However, some properties of $\mathcal{H}_1$ can be found in \cite{Zagier:75}.\\[0.2cm]
\hbox{}\quad We mention that $E'$ so defined is an Hecke eigenform with respect to the Hecke algebra $\widetilde{\mathcal{H}_v}$ for any finite $v$ not even.
This simply follows from Lemma 3.8 and the definition of $f^+_{K,v}$.\\[0.2cm]
\hbox{}\quad For any $z=(z_1,...,z_n)=(x_1+iy_1,...,x_n+iy_n)\in\mathfrak{h}^n,$ we take $g_z\in\widetilde{SL_2(\mathbb{A})}$ as $g_z=g_xg_y=\mathbf{u}^\sharp(x)\mathbf{m}(\sqrt{y})$ where
\[
x_v=\begin{cases}
0 &\mbox{ if }v<\infty,\\
x_i&\mbox{ if }v=\infty_i,
\end{cases}
\]
and
\[
y_v=\begin{cases}
1&\mbox{ if }v<\infty,\\
y_i&\mbox{ if }v=\infty_i.
\end{cases}
\]
Thus one has $g_z\mathbf{i}=z.$\\[0.2cm]
\textbf{Definition 8.3.}
For $\kappa\geq1,$ the Eisenstein series $E(z)=E_{\kappa+1/2}(z,\chi')$ on $\mathfrak{h}^n$ of weight $\kappa+1/2$ twisted by $\chi'$ is defined by
\[
E(z)=E'(g_z)\prod_{v\,|\,\infty}\tilde{j}(g_{z,v},\mathbf{i})^{2\kappa+1}
\]
where $z\in\mathfrak{h}.$\\[0.2cm]
\hbox{}\quad One can check that $E\in M_{\kappa+1/2}(\Gamma),$ hence $E'\in\mathbf{A}_{\kappa+1/2}(SL_2(F)\backslash\widetilde{SL_2(\mathbb{A})}; \widetilde{\Gamma_f},\varepsilon).$\\[0.2cm]
\hbox{}\quad In order to get the Fourier coefficients of $E(z),$ we have to calculate the Fourier coefficients of $E'$ on $g_z.$
We separate our works into two cases $\kappa\geq2$ and $\kappa=1$, though the latter case is not much different from the former.


\section{The Case $\kappa\geq 2$}
For any adele $h\in\mathbb{A}, g\in\widetilde{SL_2(\mathbb{A})},$ we have
\[
E'(\mathbf{u}^\sharp(h+x)g)=E'(\mathbf{u}^\sharp(h)g),
\]
where $x\in F$ is arbitrary chosen.
So taking $dh$ as the Haar measure of $\mathbb{A}/F$ such that Vol$(\mathbb{A}/F)=1,$ we have
\begin{equation}
\label{expansion}
E'(g_z)=E'(\mathbf{u}^\sharp(x)g_y)=\sum_{\xi\in F}c(\xi,g_y)\psi_\eta(\xi x).
\end{equation}
where
\[
c(\xi,g_y)=\int_{\mathbb{A}/F}E(\mathbf{u}^\sharp(h)g)\overline{\psi_\eta(\xi h)}dh=\sum_{\gamma\in B\backslash G}\int_{\mathbb{A}/F}f(\gamma\mathbf{u}^\sharp(h)g)\overline{\psi_\eta(\xi h)}dh
\]
is the $\xi-$th Fourier coefficient associated to $g_y.$\\
\hbox{}\quad Since $B\backslash G$ has $\{1\}\cup\{\mathbf{w}_1\mathbf{u}^\sharp(x)\,|\,x\in F\}$ as a complete system of representatives, the coefficients become
\begin{equation}
\label{coeff}
c(\xi,g_y)=\int_{\mathbb{A}/F}f(g_y)\overline{\psi_\eta(\xi h)}dh+\int_{\mathbb{A}}f(\mathbf{w}_1\mathbf{u}^\sharp(t)g_y)\overline{\psi_\eta(\xi t)}dt.
\end{equation}
Notice that the front term  of the right hand side is $f(g_y)$ or 0 if $\xi=0$ or $\xi\neq 0,$ respectively.
Let us focus on the latter term first.
Now since it is well-known that if $dt_v$ is the Haar measure for $F_v$ such that Vol($\mathfrak{o}_v$)=1 if $v<\infty$ or the usual measure for $\mathbb{R}$ if $v\,|\,\infty,$
we have
\[
dh=|\mathfrak{D}|^{-1/2}\prod_{v}dt_v
\]
where $\mathfrak{D}$ is the discriminant of $F$ over $\mathbb{Q},$ so
\[
\begin{aligned}
&\int_{\mathbb{A}}f(\mathbf{w}_1\mathbf{u}^\sharp(t)g_y)\overline{\psi_\eta(\xi t)}dt\\
=|\mathfrak{D}|^{-1/2}&\times\prod_{v<\infty}\int_{F_v}f^+_{K,v}(\mathbf{w}_1\mathbf{u}^\sharp(t))\overline{\psi_\eta(\xi t)}dt\\
&\times\prod_{j=1}^n\int_{-\infty}^{\infty}\tilde{j}(\mathbf{w}_1\mathbf{u}^\sharp(t)g_{y,\infty_j},\mathbf{i})^{-(2\kappa+1)}\mathbf{e}(-\eta\xi_{\infty_j}t)dt.
\end{aligned}
\]
If $\xi\neq 0,$ by Proposition 4.5, the finite part becomes
\[
\prod_{v<\infty}\int_{F_v}f^+_{K,v}(\mathbf{w}_1\mathbf{u}^\sharp(t))\overline{\psi_\eta(\xi t)}dt
=\prod_{v<\infty}\bgamma(\xi_v,q_v^{-s_v})q_v^{-\mathfrak{f}_{\xi_v}s_v}\Psi(\xi_v,q_v^{-s_v}).
\]
Because
\[
\bgamma(\xi_v,q_v^{-s_v})=\frac{1-q_v^{-2s_v}/q_v}{1-\chi_{\xi_v}q_v^{-s_v}/\sqrt{q_v}}=\frac{1-\chi'_v(\varpi_v)^2q_v^{-2\kappa}}{1-\chi_{\xi_v}\chi'_v(\varpi_v)q_v^{-\kappa}},
\]
we reduce the product as
\[
\prod_{v<\infty}\bgamma(\xi_v,q_v^{-s_v})^{-1}q_v^{-\mathfrak{f}_{\xi_v}s_v}\Psi(\xi_v,q_v^{-s_v})
=\frac{L_F(\kappa,\chi_\xi\chi')}{L_F(2\kappa,\chi'^2)}\prod_{v<\infty}q_v^{-\mathfrak{f}_{\xi_v}s_v}\Psi(\xi_v,q_v^{-s_v})
\]
where $L_F(\ast,\theta)$ is the L-function of $F$ associated to the character $\theta$ and $\chi_\xi$ is the character generated multiplicatively by $\chi_{\xi_v}$ defined in Section 2.
In fact, $\chi_\xi$ is the usual character of $F(\sqrt{\xi})$ over $F$.
Notice that by the definition of $\Psi$ and Lemma 7.2, we get that for nonzero $\xi,$ $c(\xi,g)=0$ unless $\xi\equiv\square$ mod 4.
We do not need to consider the case $\xi=0$ for the finite part since it is finite by Proposition 4.8 and as we will see, the infinite part for it vanishes.
\hbox{}\quad Now for the infinite part, we have
\[
\begin{aligned}
&\int_{-\infty}^{\infty}\tilde{j}(\mathbf{w}_1\mathbf{u}^\sharp(t)g_{y,\infty_j},\mathbf{i})^{-(2\kappa+1)}\mathbf{e}(-\eta\xi_{\infty_j}t)dt\\
=&\int_{-\infty}^{\infty}(i\sqrt{y_j}+t/\sqrt{y_j})^{-(\kappa+1/2)}e^{-2\pi i\eta\xi_{\infty_j}t}dt\\
=&i^{-(\kappa+3/2)}y_j^{\kappa/2+1/4}e^{-2\pi\eta\xi_{\infty_j}y_j}\int_{y_j-i\infty}^{y_j+i\infty}t^{-(\kappa+1/2)}e^{2\pi\eta\xi_{\infty_j}t}dt\\
=&\begin{cases}
(-2\pi i)^{\kappa+1/2}y_j^{\kappa/2+1/4}e^{-2\pi\eta\xi_{\infty_j}y_j}(\eta\xi_{\infty_j})^{\kappa-1/2}\Gamma(\kappa+1/2)^{-1}&\mbox{ if }\eta\xi_{\infty_j}>0\\
0&\mbox{ if }\eta\xi_{\infty_j}\leq 0\\
\end{cases}
\end{aligned}
\]
by the inverse Laplace transform of $t^{-(\kappa+1/2)}.$
Thus if $\eta\xi\succ 0,$ i.e., if $\eta\xi$ is totally positive, the infinite part becomes
\[
\begin{aligned}
&\prod_{j=1}^n\int_{-\infty}^{\infty}\tilde{j}(\mathbf{w}_1\mathbf{u}^\sharp(t)g_{y,\infty_j},\mathbf{i})^{-(2\kappa+1)}\mathbf{e}(-\eta\xi_{\infty_j}t)dt\\
=&(-2\pi i)^{n(\kappa+1/2)}N(y_\infty)^{\kappa/2+1/4}\mathbf{e}(\eta\xi\ast iy_\infty)N_{F/\mathbb{Q}}(\eta\xi)^{\kappa-1/2}\Gamma(\kappa+1/2)^{-n}.
\end{aligned}
\]
Here $y_\infty$ denotes the infinite part of $y$ and $N(y_\infty)=\prod_{j=1}^n y_j.$
Now for nonzero $\xi,$ we have
\begin{equation}
\label{equa}
\begin{aligned}
&c(\xi,g_y)\\
=&|\mathfrak{D}|^{-1/2}\frac{L_F(\kappa,\chi_\xi\chi')}{L_F(2\kappa,\chi'^2)}\prod_{v<\infty}q_v^{-\mathfrak{f}_{\xi_v}s_v}\Psi(\xi_v,q_v^{-s_v})\\
&\times(-2\pi i)^{n(\kappa+1/2)}N(y_\infty)^{\kappa/2+1/4}\mathbf{e}(\eta\xi\ast iy_\infty)N_{F/\mathbb{Q}}(\eta\xi)^{\kappa-1/2}\Gamma(\kappa+1/2)^{-n}
\end{aligned}
\end{equation}
and this only occurs when $\xi\equiv\square$ mod 4 and $\eta\xi$ is totally positive.\\
\hbox{}\quad For $\xi=0,$ we have just seen that the latter term in equation (\ref{coeff}) vanishes, thus
\[
c(0,g_y)=f(g_y).
\]
By an easy calculation, we get that for $v<\infty$
\[
f_{K,v}^+(1)=q_v^{(2s_v+1/2)e_v+c_vs_v}\alpha_\psi(\eta)=\chi'_v(4\bdelta_v)^{-1}|2|_v^{1/2-2\kappa}|\bdelta_v|^{1/2-\kappa}\alpha_{\psi_v}(\eta).
\]
So
\[
c(0,g_y)=N(y_\infty)^{\kappa/2+1/4}\prod_{v<\infty}\chi'_v(4\bdelta_v)^{-1}|2|_v^{1/2-2\kappa}|\bdelta_v|^{1/2-\kappa}\alpha_{\psi_v}(\eta).
\]
Notice that 
\[
\prod_{v<\infty}\chi'_v(4\bdelta_v)^{-1}|2|_v^{1/2-2\kappa}|\bdelta_v|^{1/2-\kappa}=\chi'(\mathfrak{d})^{-1}2^{n(2\kappa-1/2)}|\mathfrak{D}|^{\kappa-1/2}
\]
and
\[
\prod_{v<\infty}\alpha_\psi(\eta)=\prod_{v|\infty}\alpha_\psi(-\eta)=\prod_{v|\infty}\alpha_\psi((-1)^{\kappa+1})=\exp\left(\frac{(-1)^{\kappa+1}n\pi\sqrt{-1}}{4}\right),
\]
we can rewrite $c(0,g_y)$ as
\[
c(0,g_y)=2^{n(2\kappa-1/2)}|\mathfrak{D}|^{\kappa-1/2}N(y_\infty)^{\kappa/2+1/4}\chi'(\mathfrak{d})^{-1}\exp\left(\frac{(-1)^{\kappa+1}n\pi\sqrt{-1}}{4}\right).
\]
\hbox{}\quad Substituting our results into equation (\ref{expansion}), we get
\[
\begin{aligned}
E'(g_z)=&2^{n(2\kappa-1/2)}|\mathfrak{D}|^{\kappa-1/2}N(y_\infty)^{\kappa/2+1/4}\chi'(\mathfrak{d})^{-1}\exp\left(\frac{(-1)^{\kappa+1}n\pi\sqrt{-1}}{4}\right)\\
&+\sum_{\begin{smallmatrix}\xi\equiv\square\,\mathrm{mod}\,4\\ \eta\xi\succ 0\end{smallmatrix}}|\mathfrak{D}|^{-1/2}\frac{L_F(\kappa,\chi_\xi\chi')}{L_F(2\kappa,\chi'^2)}\prod_{v<\infty}q_v^{-\mathfrak{f}_{\xi_v}s_v}\Psi(\xi_v,q_v^{-s_v})\\
&\times(-2\pi i)^{n(\kappa+1/2)}N(y_\infty)^{\kappa/2+1/4}q^{\eta\xi}N_{F/\mathbb{Q}}(\eta\xi)^{\kappa-1/2}\Gamma(\kappa+1/2)^{-n}\\[0.2cm]
=&2^{n(2\kappa-1/2)}|\mathfrak{D}|^{\kappa-1/2}N(y_\infty)^{\kappa/2+1/4}\chi'(\mathfrak{d})^{-1}\exp\left(\frac{(-1)^{\kappa+1}n\pi\sqrt{-1}}{4}\right)\\
&+\sum_{\begin{smallmatrix}\eta\xi\equiv\square\,\mathrm{mod}\,4\\ \xi\succ 0\end{smallmatrix}}|\mathfrak{D}|^{-1/2}\frac{L_F(\kappa,\chi_{\eta\xi}\chi')}{L_F(2\kappa,\chi'^2)}\prod_{v<\infty}q_v^{-\mathfrak{f}_{\eta\xi_v}s_v}\Psi(\eta\xi_v,q_v^{-s_v})\\
&\times(-2\pi i)^{n(\kappa+1/2)}N(y_\infty)^{\kappa/2+1/4}q^\xi N_{F/\mathbb{Q}}(\xi)^{\kappa-1/2}\Gamma(\kappa+1/2)^{-n}
\end{aligned}
\]
Now we can get the Fourier coefficients of $E(z).$
Since
\[
E(z)=E'(g_z)\prod_{v\,|\,\infty}\tilde{j}(g_{z,v},\mathbf{i})^{2\kappa+1}=E'(g_z)N(y_\infty)^{-(\kappa/2+1/4)},
\]
we have
\[
\begin{aligned}
E(z)=&2^{n(2\kappa-1/2)}|\mathfrak{D}|^{\kappa-1/2}\chi'(\mathfrak{d})^{-1}\exp\left(\frac{(-1)^{\kappa+1}n\pi\sqrt{-1}}{4}\right)\\
&+\sum_{\begin{smallmatrix}\eta\xi\equiv\square\,\mathrm{mod}\,4\\ \xi\succ 0\end{smallmatrix}}|\mathfrak{D}|^{-1/2}\frac{L_F(\kappa,\chi_{\eta\xi}\chi')}{L_F(2\kappa,\chi'^2)}\prod_{v<\infty}q_v^{-\mathfrak{f}_{\eta\xi_v}s_v}\Psi(\eta\xi_v,q_v^{-s_v})\\
&\times(-2\pi i)^{n(\kappa+1/2)}N_{F/\mathbb{Q}}(\xi)^{\kappa-1/2}\Gamma(\kappa+1/2)^{-n}q^\xi.
\end{aligned}
\]
From this we see that $E(z)$ actually lies in $M^+_{\kappa+1/2}(\Gamma).$ But the formula is still a little messy. However, by the functional equation of L-functions of $F$ and $F(\sqrt{\eta\xi})$
and the identity 
\[L_F(w,\chi_{\eta\xi}\chi')L_F(w,\chi')=L_{F(\sqrt{\eta\xi})}(w,\chi'\circ\mbox{Norm}_{F(\sqrt{\eta\xi})/F})
\]
from class field theory, we have
\[
L_F(2\kappa,\chi'^2)=|\mathfrak{D}|^{1/2-2\kappa}\chi'(\mathfrak{d})^2\pi^{n(2\kappa-1/2)}\left(\frac{\Gamma(1/2-\kappa)}{\Gamma(\kappa)}\right)^nL_F(1-2\kappa,\overline{\chi'}^2)
\]
and
\[
\begin{aligned}
&L_F(\kappa,\chi_{\eta\xi}\chi')\\
=&\begin{cases}
|\mathfrak{D}N_{F/\mathbb{Q}}(\mathfrak{D}_{\eta\xi})|^{1/2-\kappa}\chi'(\mathfrak{D}_{\eta\xi}\mathfrak{d})\pi^{n(\kappa-1/2)}\left(\frac{\Gamma((1-\kappa)/2)}{\Gamma(\kappa/2)}\right)^nL_F(1-\kappa,\overline{\chi_{\eta\xi}\chi'})&\mbox{ if $\kappa$ is even,}\\
|\mathfrak{D}N_{F/\mathbb{Q}}(\mathfrak{D}_{\eta\xi})|^{1/2-\kappa}\chi'(\mathfrak{D}_{\eta\xi}\mathfrak{d})(4\pi)^{n(\kappa-1/2)}\left(\frac{\Gamma(1-\kappa)\Gamma(\kappa/2)}{\Gamma(\kappa)\Gamma((1-\kappa)/2)}\right)^nL_F(1-\kappa,\overline{\chi_{\eta\xi}\chi'})&\mbox{ if $\kappa$ is odd,}\\
\end{cases}
\end{aligned}
\]
where $\mathfrak{D}_{\eta\xi}$ is the relative discriminant of $F(\sqrt{\eta\xi})$ over $F.$
Substituting these identities back to the Fourier expansion of $E(z)$ and using the properties of the gamma function, we get
\[
\begin{aligned}
E(z)=&2^{n(2\kappa-1/2)}|\mathfrak{D}|^{\kappa-1/2}\chi'(\mathfrak{d})^{-1}\exp\left(\frac{(-1)^{\kappa+1}n\pi\sqrt{-1}}{4}\right)\\
&+\sum_{\begin{smallmatrix}\eta\xi\equiv\square\,\mathrm{mod}\,4\\ \xi\succ 0\end{smallmatrix}}|\mathfrak{D}|^{\kappa-1/2}\chi'(\mathfrak{D}_{\eta\xi}\mathfrak{d}^{-1})\frac{L_F(1-\kappa,\overline{\chi_{\eta\xi}\chi'})}{L_F(1-2\kappa,\overline{\chi'}^2)}\prod_{v<\infty}q_v^{-\mathfrak{f}_{\eta\xi_v}s_v}\Psi(\eta\xi_v,q_v^{-s_v})\\
&\times\exp\left(\frac{(-1)^{\kappa+1}n\pi\sqrt{-1}}{4}\right)2^{n(2\kappa-1/2)}\left|\frac{N_{F/\mathbb{Q}}(\xi)}
{N_{F/\mathbb{Q}}(\mathfrak{D_{\eta\xi}})}\right|^{\kappa-1/2}q^\xi.
\end{aligned}
\]
\quad Finally, let $G(z)$ be the Eisenstein series normalized from $E(z)$ such that the constant term of whose Fourier series is $L_F(1-2\kappa,\overline{\chi'}^2),$ then the explicit form of $G(z)$ is as the following theorem.\\[0.2cm]
\textbf{Theorem 9.1.}
Let $\eta=(-1)^\kappa.$
For $\kappa\geq2,$ the Fourier series expansion of $G(z)=G_{\kappa+1/2}(z,\chi')$ defined above is
\[
G(z)=L_F(1-2\kappa,\overline{\chi'}^2)+\sum_{\begin{smallmatrix}(-1)^{\kappa}\xi\equiv\square\,\mathrm{mod}\,4\\ \xi\succ 0\end{smallmatrix}}\chi'(\mathfrak{D}_{\eta\xi})L_F(1-\kappa,\overline{\chi_{\eta\xi}\chi'})\mathfrak{C}_\kappa(\eta\xi)q^\xi,
\]
where
\[
\begin{aligned}
\mathfrak{C}_\kappa(\xi)
=&N_{F/\mathbb{Q}}(\mathfrak{F}_\xi)^{\kappa-1/2}\chi'(\mathfrak{F}_\xi)\prod_{v<\infty}\Psi(\xi_v,q_v^{1/2-\kappa}\chi'_v(\varpi_v))\\
=&\sum_{\mathfrak{a}|\mathfrak{F}_\xi}\mu(\mathfrak{a})\chi_\xi(\mathfrak{a})\chi'(\mathfrak{a})N_{F/\mathbb{Q}}(\mathfrak{a})^{\kappa-1}\sigma_{2\kappa-1,\chi'^2}(\mathfrak{F}_\xi\mathfrak{a}^{-1}).
\end{aligned}
\]
Here $\mathfrak{F}_\xi^2\mathfrak{D_\xi}=(\xi),$ $\mathfrak{a}$ runs over all integral ideals dividing $\mathfrak{F_\xi},$ $\mu$ is the M\"{o}bius function and $\sigma_{k,\chi}$ is the sum of divisors function twisted by $\chi,$ that is,
\[
\sigma_{k,\chi}(\mathfrak{A})=\sum_{\mathfrak{b}|\mathfrak{A}}N_{F/\mathbb{Q}}(\mathfrak{b})^k\chi(\mathfrak{b})
\]
for any integral ideal $\mathfrak{A}$ of $F.$
Moreover, $G$ is a Hecke eigenform with respect to the Hecke algebra $\widetilde{\mathcal{H}_v}$ for any finite $v$ which is not even.\\[0.2cm]\hbox{}\quad One can easily see that if $F=\mathbb{Q},$ the Fourier series $G_{\kappa+1/2}$ coincides with $\mathcal{H}_\kappa$ introduced in \cite{Cohen:75}.


\section{The Case $\kappa=1$}
In this section $F\neq\mathbb{Q}.$
We use the same notations as in Section 8.
Similarly, for $z\in\mathfrak{h}^n$, we have
\begin{equation}
\label{expansion2}
E'_\epsilon(g_z)=E'_\epsilon(\mathbf{u}^\sharp(x)g_y)=\sum_{\xi\in F}c_\epsilon(\xi,g_y)\psi_\eta(\xi x).
\end{equation}
where
\begin{equation}
\label{coeff2}
c_\epsilon(\xi,g_y)=\int_{\mathbb{A}/F}f_\epsilon(g_y)\overline{\psi_\eta(\xi h)}dh+\int_{\mathbb{A}}f_\epsilon(\mathbf{w}_1\mathbf{u}^\sharp(t)g_y)\overline{\psi_\eta(\xi t)}dt.
\end{equation}
For simplicity, denote the former of the right hand by $A(\xi,\epsilon)$ and the latter by $B(\xi,\epsilon).$
What we want to obtain are the values of their analytic continuation at $\epsilon=0$.
Let us focus on $B(\xi,\epsilon)$ first.
We have
\[
\begin{aligned}
&B(\xi,\epsilon)\\
=&|\mathfrak{D}|^{-1/2}\times\prod_{v<\infty}\int_{F_v}f^+_{K,v}(\mathbf{w}_1\mathbf{u}^\sharp(t))\overline{\psi_\eta(\xi t)}dt\\
&\times\prod_{j=1}^n\int_{-\infty}^{\infty}\tilde{j}(\mathbf{w}_1\mathbf{u}^\sharp(t)g_{y,\infty_j},\mathbf{i})^{-3}|\tilde{j}(\mathbf{w}_1\mathbf{u}^\sharp(t)g_{y,\infty_j},\mathbf{i})|^{-2\epsilon}\mathbf{e}(-\eta\xi_{\infty_j}t)dt.
\end{aligned}
\]
\hbox{}\quad If $\xi\neq0,$ the calculation of the finite part is nothing different from which in Section 9 except for $s_v$ being changed into $s_v+\epsilon,$ thus
\[
\prod_{v<\infty}\int_{F_v}f^+_{K,v}(\mathbf{w}_1\mathbf{u}^\sharp(t))\overline{\psi_\eta(\xi t)}dt
=\frac{L_F(1+\epsilon,\chi_\xi\chi')}{L_F(2+2\epsilon,\chi'^2)}\prod_{v<\infty}q_v^{-\mathfrak{f}_{\xi_v}s_v+\epsilon}\Psi(\xi_v,q_v^{-s_v-\epsilon}).
\]
This is not zero only if $\xi\equiv\square$ mod 4 and has a simple pole at $\epsilon=0$ if $\chi_\xi\chi'=1.$
And this happen only if $\chi_\xi$ is unramified at infinite places.\\
\hbox{}\quad If $\xi=0,$ we have
\[
\prod_{v<\infty}\int_{F_v}f^+_{K,v}(\mathbf{w}_1\mathbf{u}^\sharp(t))dt=\frac{L_F(1+2\epsilon,\chi'^2)}{L_F(2+2\epsilon,\chi'^2)}.
\]
from Proposition 4.8.
Again, this has a simple pole at $\epsilon=0$ if $\chi'2=1.$\\
\hbox{}\quad To calculate the infinite part, we need some basic properties of the confluent hypergeometric functions which will be used.\\[0.2cm]
\textbf{Lemma 10.1.}
Let $\alpha, \beta\in\mathbb{C}$ and $h\in\mathbb{R}$ such that Re$(\alpha)>0,$ Re$(\beta)>0$ and Re$(\alpha+\beta)>1.$
Put $u=t+is\in\mathfrak{h}$ with a fixed $s>0.$
Then we have
\[
\begin{aligned}
&\int^\infty_{-\infty}u^{-\alpha}\bar{u}^{-\beta}e^{-2\pi iht}dt\\
=&(2\pi)^{\alpha+\beta}i^{\beta-\alpha}\Gamma(\alpha)^{-1}\Gamma(\beta)^{-1}e^{2\pi hs}\times
\begin{cases}
h^{\alpha+\beta-1}e^{-4\pi hs}\sigma(4\pi hs,\alpha,\beta)&\mbox{ if }h>0,\\
\Gamma(\alpha+\beta-1)(4\pi s)^{1-\alpha-\beta}&\mbox{ if }h=0,\\
|h|^{\alpha+\beta-1}\sigma(4\pi|h|s,\beta,\alpha)&\mbox{ if }h<0,\\
\end{cases}
\end{aligned}
\]
where
\begin{equation}
\label{sigma}
\sigma(s,\alpha,\beta)=\int^\infty_0(u+1)^{\alpha-1}u^{\beta-1}e^{-su}du
\end{equation}
for $s>0.$\\
\textit{Proof.} See the proof of Lemma 1 in \cite{Shimura:75}.\done\\[0.2cm]
\textbf{Lemma 10.2.} If we put
\[
V(s,\alpha,\beta)=e^{-s/2}\Gamma(\beta)^{-1}s^\beta\sigma(s,\alpha,\beta),
\]
then $V(s,\alpha,\beta)$ can be continued as a holomorphic function in $(\alpha,\beta)$ to the whole $\mathbb{C}^2$ and satisfies the functional equation
\[
V(s,1-\beta,1-\alpha)=V(s,\alpha,\beta).
\]
Thus $\sigma(s,\alpha,\beta)$ also can be continued as a holomorphic function in $(\alpha,\beta)$ to the whole $\mathbb{C}^2$.
Moreover, we have
\[
\sigma(s,1,\beta)=\Gamma(\beta)s^{-\beta}.
\]
\textit{Proof.} See A3 of \cite{Shimura:07}.\done\\[0.2cm]
\hbox{}\quad Now for an infinite place $\infty_j,$ we have
\[
\begin{aligned}
&\int_{-\infty}^{\infty}\tilde{j}(\mathbf{w}_1\mathbf{u}^\sharp(t)g_{y,\infty_j},\mathbf{i})^{-3}|\tilde{j}(\mathbf{w}_1\mathbf{u}^\sharp(t)g_{y,\infty_j},\mathbf{i})|^{-2\epsilon}\mathbf{e}(-\eta\xi_{\infty_j}t)dt\\
=&\int_{-\infty}^{\infty}(i\sqrt{y_j}+t/\sqrt{y_j})^{-3/2}|i\sqrt{y_j}+t/\sqrt{y_j}|^{-\epsilon}e^{-2\pi i\eta\xi_{\infty_j}t}dt\\
=&y_j^{3/4+\epsilon}\int_{-\infty}^{\infty}(iy_j+t)^{-3/2-\epsilon/2}(-iy_j+t)^{-\epsilon/2}e^{-2\pi i\eta\xi_{\infty_j}t}dt\\
=&(2\pi)^{3/2+\epsilon}i^{-3/2}\Gamma(3/2+\epsilon/2)^{-1}\Gamma(\epsilon/2)^{-1}e^{-2\pi\xi_{\infty_j}y_j}y_j^{3/4+\epsilon/2}\\
&\times\begin{cases}
|\xi_{\infty_j}|^{1/2+\epsilon}e^{4\pi i\eta\xi_{\infty_j}t}\sigma(-4\pi\xi_{\infty_j}y_j,3/2+\epsilon/2,\epsilon/2)&\mbox{ if }\xi_{\infty_j}<0,\\
\Gamma(1/2+\epsilon)(4\pi y_j)^{-1/2-2\epsilon/2}&\mbox{ if }\xi_{\infty_j}=0,\\
\xi_{\infty_j}^{1/2+\epsilon}\sigma(4\pi\xi_{\infty_j}y_j,\epsilon/2,3/2+\epsilon/2)&\mbox{ if }\xi_{\infty_j}>0\\
\end{cases}
\end{aligned}
\]
by Lemma 10.1.
Notice that here $\eta=-1.$\\
\hbox{}\quad Using the fact that the gamma function has a pole at $0$ and the integral in equation (\ref{sigma}) defining $\sigma(s,\alpha,\beta)$ converges for Re$(\beta)>0,$
we see the infinite part vanishes at $\epsilon=0$ if $\xi_{\infty_j}\geq 0.$
Note that if the finite part we just calculated has a simple pole, then $\xi=0$ or $\chi_\xi$ is unramified at all infinite places, that is, $\xi\succ0.$
In both of the cases, the infinite part gives a zero of order $n\geq2.$
So we see that $B(\xi,0)$ only occurs for $-\xi\succ0.$\\
\hbox{}\quad If $\xi_{\infty_j}<0$, by Lemma 10.2,
\[
\begin{aligned}
&\sigma(-4\pi\xi_{\infty_j}y_j,3/2+\epsilon/2,\epsilon/2)\\
=&\Gamma(\epsilon/2)\Gamma(-1/2-\epsilon/2)^{-1}|4\pi\xi_{\infty_j}y_j|^{-1/2-3\epsilon/2}\sigma(-4\pi\xi_{\infty_j}y_j,1-\epsilon/2,-1/2-\epsilon),
\end{aligned}
\]
so the integral for $\infty_j$ has the limit
$(-2\pi i)^{3/2}\Gamma(3/2)^{-1}e^{2\pi\xi_{\infty_j}y_j}y_j^{3/4}|\xi_{\infty_j}|^{1/2}$ at $\epsilon=0$.\\
\hbox{}\quad Substituting these results back to $B(\xi,\epsilon)$ and let $\epsilon$ approach 0,
we get
\[
\begin{aligned}
&B(\xi,0)\\
=&\begin{cases}
|\mathfrak{D}|^{-1/2}(-2\pi i)^{3n/2}N(y_\infty)^{3/4}\mathbf{e}(-\xi\ast iy_\infty)N_{F/\mathbb{Q}}(-\xi)^{1/2}\Gamma(3/2)^{-n}\\
\times\frac{L_F(1,\chi_\xi\chi')}{L_F(2,\chi'^2)}\prod_{v<\infty}q_v^{-\mathfrak{f}_{\xi_v}s_v}\Psi(\xi_v,q_v^{-s_v})\\
\,\,\,\,\mbox{ if }\xi\equiv\square\mbox{ mod 4 and }-\xi\succ0,\\
0\mbox{ otherwise,}
\end{cases}
\end{aligned}
\]
which coincides with equation (\ref{equa}).\\
\hbox{}\quad The calculation of $A(\xi,\epsilon)$ is similar to which of the front term in the right hand side of equation (\ref{coeff}).
We get
\[
A(\xi,0)=
\begin{cases}
2^{3n/2}|\mathfrak{D}|^{1/2}N(y_\infty)^{3/4}\chi'(\mathfrak{d})^{-1}\exp(n\pi i/4)&\mbox{ if }\xi=0,\\
0&\mbox{ if }\xi\neq0.
\end{cases}
\]
Then by exactly the same argument as in the last section, we can get the similar theorem with Theorem 9.1 for the case $\kappa=1.$\\[0.2cm]
\textbf{Theorem 10.3.}
Let $G(z)=G_{3/2}(z,\chi')$ be the Eisenstein normalized from $E(z)=E_{3/2}(z,\chi')$ such that the constant term of whose Fourier series is $L_F(-1,\overline{\chi'}^2).$
We have
\[
G(z)=L_F(-1,\overline{\chi'}^2)+\sum_{\begin{smallmatrix}-\xi\equiv\square\,\mathrm{mod}\,4\\ \xi\succ 0\end{smallmatrix}}\chi'(\mathfrak{D}_{\eta\xi})L_F(0,\overline{\chi_{-\xi}\chi'})\mathfrak{C}_1(-\xi)q^\xi,
\]
where
\[
\begin{aligned}
\mathfrak{C}_1(\xi)
=&N_{F/\mathbb{Q}}(\mathfrak{F}_\xi)^{1/2}\chi'(\mathfrak{F}_\xi)\prod_{v<\infty}\Psi(\xi_v,q_v^{-1/2}\chi'_v(\varpi_v))\\
=&\sum_{\mathfrak{a}|\mathfrak{F}_\xi}\mu(\mathfrak{a})\chi_\xi(\mathfrak{a}){\chi'(\mathfrak{a})}\sigma_{1,{\chi'^2}}(\mathfrak{F}_\xi\mathfrak{a}^{-1}).
\end{aligned}
\]
Moreover, $G$ is a Hecke eigenform with respect to the Hecke algebra $\widetilde{\mathcal{H}_v}$ for any finite $v$ which is not even.\\[0.2cm]
\textit{Remark.}
If $F=\mathbb{Q}$, for any integer $m\in Z$, the finite part of $B(m^2,\epsilon)$ gives a simple pole while the infinite part only gives a simple zero at $\epsilon=0.$
Thus $B(m^2,0)$ does not vanish and our calculation does not work any more.


\section{A Corollary}
In this section, we consider the same condition as in Section 9 and let $\kappa$ be a positive integer.
In addition, we set the three groups
\[
U=\left\{\begin{pmatrix}1&x\\0&1\end{pmatrix}\in SL_2(\mathbb{A})\,\bigg|\, x\in\mathbb{A}\right\},
\]
\[
D=\left\{\begin{pmatrix}a&0\\0&a^{-1}\end{pmatrix}\in SL_2(\mathbb{F})\,\bigg|\, a\in F^\times\right\}
\]
and
\[
\Xi=\prod'_{v<\infty}\Gamma[(4\mathfrak{d_v})^{-1},4\mathfrak{d_v}].
\]
Here the restricted product $\prod'$ is respected to $SL_2(\mathfrak{o_v}).$\\[0.2cm]
\textbf{Lemma 11.1.}
Let $f\in M^+_{\kappa+1/2}(\Gamma)$ and $\varphi=\varphi_f$ be its corresponding automorphic form in $\mathbf{A}_{\kappa+1/2}(SL_2(F)\backslash\widetilde{SL_2(\mathbb{A})}; \widetilde{\Gamma_f},\varepsilon),$
then the constant term in the Fourier expansion of $\varphi$ is determined by its values on
\[
UD\backslash\widetilde{SL_2(\mathbb{A})}/\Xi\prod_{j=1}^n\widetilde{SL_2(F_{\infty_j})}.
\]
\textit{Proof.}
Denote the constant term by $\varphi_0$ and let $g\in\widetilde{SL_2(\mathbb{A})}.$
Then
\[
\varphi_0(g)=\int_{\mathbb{A}/F}\varphi(\mathbf{u}^\sharp(x)g)dx
\]
for $g\in\widetilde{SL_2(\mathbb{A})}.$
It is trivial that $\varphi_0(\mathbf{u}^\sharp(x)\mathbf{m}(a)g)=\varphi_0(g)$ for $x\in\mathbb{A}, a\in F^\times.$
Now we decompose $g$ into $g=g_f\prod_{j=1}^ng_{\infty_j}$ such that $g_{\infty_j}$ is the $\infty_j$-th component of $g$ and $g_{\infty_j}=[g_{\infty_j},1]\in\widetilde{SL_2(F_{\infty_j})}=\widetilde{SL_2(\mathbb{R})}.$
By Iwasawa decomposition, $g_{\infty_j}$ has the form
\[
g_{\infty_j}=\begin{pmatrix}a_j&b_j\\0&a_j^{-1}\end{pmatrix}\begin{pmatrix}\cos\theta_j&-\sin\theta_j\\\sin\theta_j&\cos\theta_j\end{pmatrix}
\]
for some $a_j, b_j, \theta_j\in\mathbb{R}.$
From Lemma 6.2 we get
\[
\varphi_0(g)=\varphi_0(g_f)\prod_{j=1}^{n}\frac{\alpha_{\psi_{\infty_j}}(\eta)}{\alpha_{\psi_{\infty_j}}(\eta a_j)}|a_j|^{\kappa+1/2}e^{2\pi i(\kappa+1/2)\theta_j}.
\]
So below we can only consider the case $g_\infty=1.$
By Theorem 7.4, we know $\rho(E^K)\varphi=\varphi.$
From the integral formula of $\varphi_0$ one easily get that also $\rho(E^K)\varphi_0=\varphi_0.$
To show that $\varphi_0$ only depends on its values on $\widetilde{SL_2(\mathbb{A}_f)}/\Xi,$ it suffices to show that if $m$ is a representative of $\widetilde{SL_2(\mathbb{A}_f)}/\Xi$,
the space of $\varphi_0(m\gamma)$ considered as a function in $\gamma\in\Xi$ where $\varphi_0$ is any constant term we can take in this proof is of only one dimension.
By Iwasawa decomposition again, we can always take $m$ to be upper triangular.
Hence we only have to show that for $v<\infty,$ the space
\[
S=\{h:\mbox{genius function on }\widetilde{\Gamma[(4\mathfrak{d_v})^{-1},4\mathfrak{d_v}]}\,|\,\rho(E^K)h=h,\,h(\mathbf{u}^\sharp(x)g)=h(g)\}
\]
is of only one dimension over $\mathbb{C}.$
Let $h\in S.$
For simplicity, we omit the lower subscript $v.$
Let $\widehat{\Omega}_\psi$ be the representation of $\widetilde{\Gamma[(4\mathfrak{d})^{-1},4\mathfrak{d}]}=\mathbf{w}_{2\bdelta}^{-1}\widetilde{\Lambda}\mathbf{w}_{2\bdelta}$ on $\omega_\psi(\mathbf{w}_{2\bdelta}^{-1})\mathcal{S}(\mathfrak{p}^{-e}/\mathfrak{o})$
as a conjugate of $\Omega_\psi$ defined in Section 6.
For any $\phi_\lambda$ defined by $\phi_\lambda(t)=\phi_0(t-\lambda),$ we have
\[
\omega_\psi(\mathbf{w}_{2\bdelta}^{-1})\phi_\lambda(t)=\alpha_\psi(2\bdelta)\psi(t\lambda/\bdelta)\phi_0(t),
\]
so
\[
\omega_\psi(\mathbf{w}_{2\bdelta}^{-1})\mathcal{S}(\mathfrak{p}^{-e}/\mathfrak{o})=\bigoplus_{\lambda\in\mathfrak{p}^{-e}/\mathfrak{o}}\mathbb{C}\cdot\phi'_\lambda.
\]
where $\phi'_\lambda(t)=\psi(t\lambda/\bdelta)\phi_0(t).$
Now the identity $h\ast\overline{E^K}=\rho(E^K)h=h$ tells us that $h$ is a matrix coefficient on $\Gamma[(4\mathfrak{d})^{-1},4\mathfrak{d}].$
So $h$ has the form
\[
h(g)=\sum_{\lambda_1,\lambda_2}c_{\lambda_1,\lambda_2}\overline{(\phi'_{\lambda_1},\widehat{\Omega}_\psi(g)\phi'_{\lambda_2})}
\]
where $\lambda_1, \lambda_2$ run over a complete system of representatives for $\mathfrak{p}^{-e}/\mathfrak{o}$ and $c_{\lambda_1,\lambda_2}\in\mathbb{C}.$
Then
\[
h(gk)=\sum_{\lambda_1,\lambda_2}c_{\lambda_1,\lambda_2}\left(\sum_{\lambda_3}\overline{(\phi'_{\lambda_1},\widehat{\Omega}_\psi(g)\phi'_{\lambda_3})}\overline{(\phi'_{\lambda_3},\widehat{\Omega}_\psi(k)\phi'_{\lambda_2})}\right).
\]
Using the identity $h\ast\overline{E^K}=h$ again and by Schur orthogonal relation, we can reduce $h$ to
\[
h(g)=\sum_\lambda c_\lambda\overline{(\phi'_\lambda,\widehat{\Omega}_\psi(g)\phi'_0)}.
\]
where $c_\lambda=c_{\lambda,0}.$
Let $\phi'=\sum_\lambda\overline{c_\lambda}\phi'_\lambda.$ So $h(g)=\overline{(\phi',\widehat{\Omega}_\psi(g)\phi'_0)}.$
Now for each $x\in\mathfrak{o},$ from the condition that $h(\mathbf{u}^\sharp(x/4\bdelta)g)=h(g)$ for any $g$ and the irreducibility of $\widehat{\Omega}_\psi,$ we have
\[
\psi(-xt^2/4\bdelta)\phi'(t)=\widehat{\Omega}_\psi(\mathbf{u}^\sharp(-x/4\bdelta))\phi'(t)=\phi'(t)
\]
for any $t\in\mathfrak{p}^{-e}.$
But if $\ord(t)<e,$ there exists some $x$ such that $\psi(-xt^2/4\bdelta)\neq0.$
Thus $\phi'(t)=0$ for $\ord(t)<e.$
Now by the definition of $\phi'_\lambda$ for any $\lambda\in\mathfrak{p}^{-e}/\mathfrak{o},$ we get that $\phi'$ is some constant on $\mathfrak{p}^e.$
This shows that the space $S$ of $h$ is of one dimension.
And the lemma follows.\done\\[0.2cm]
\textbf{Lemma 11.2.}
The order of the double coset space in the statement of Lemma 11.1 is the class number $h$ of $F.$\\
\textit{Proof.}
By Iwasawa decomposition, we have
\[
SL_2(\mathbb{A})=B(\mathbb{A})\Xi\prod_{\j=1}^nSL_2(F_{\infty_j})
\]
where $B(A)$ is the group of upper triangular matrices in $SL_2(\mathbb{A})$.
So one easily get that
\[
UD\backslash\widetilde{SL_2(\mathbb{A})}/\Xi\prod_{j=1}^n\widetilde{SL_2(F_{\infty_j})}\,\cong\,UD\backslash\widetilde{B(\mathbb{A})}/\widetilde{B(\mathbb{A})}\cap(\Xi\prod_{j=1}^n\widetilde{SL_2(F_{\infty_j})}).
\]
Hence it suffices to calculate the order of the double coset space in the right hand side.
Now choose a complete system $\{\mathfrak{a}_1,...,\mathfrak{a}_h\}$ of representatives for the ideal classes of $F$ and
set $a_1,...,a_h\in\mathbb{A}^\times$ be $h$ ideles such that 
\[
\begin{aligned}
\ord(a_{j,v})&=\ord(\mathfrak{a}_{j,v})&\mbox{ if }v<\infty,\\
a_{j,v}&=1&\mbox{ if }v\,|\,\infty.
\end{aligned}
\]
We show that $\mathbf{m}(a_1),...,\mathbf{m}(a_h)$ forms a complete system of representatives for the double cosets.
Obviously they are in mutually distinct double cosets.
Now for any $\beta\in\widetilde{B(\mathbb{A})},$ $\beta$ is equivalent to some $\mathbf{m}(b)$ such that $b_\infty=1.$
Let $\mathfrak{b}$ be the fractional ideal defined by $\ord(\mathfrak{b}_v)=\ord(b_v)$ for each $v<\infty.$ 
Then $\mathfrak{b}=(\gamma)\mathfrak{a}_j$ for some $\gamma\in F, 1\leq j\leq h.$
So for each $v<\infty,$ we can find some $\epsilon_v\in\mathfrak{o}^\times_v$ such that $\gamma_v^{-1}b_v\epsilon_v=a_{j,v}.$
Now it follows that $\mathbf{m}(b)$ is equivalent to $\mathbf{m}(a_j).$
So $\beta$ is equivalent to $\mathbf{m}(a_j)$ and thus the lemma.\done\\[0.2cm]
\hbox{}\quad Note that $\mathbf{m}(a_1),...,\mathbf{m}(a_h)$ used in the proof of Lemma 11.2 also forms a complete system of representatives of $UD\backslash\widetilde{SL_2(\mathbb{A})}/\Xi\prod_{j=1}^n\widetilde{SL_2(F_{\infty_j})}.$\\[0.2cm]
\textbf{Corollary 11.3.}
The Kohnen plus space $M^+_{\kappa+1/2}(\Gamma)$ is a vector space over $\mathbb{C}$ spanned by cusp forms and the $h$ Eisenstein series we constructed in Section 8, that is,
\[
M^+_{\kappa+1/2}(\Gamma)=S^+_{\kappa+1/2}(\Gamma)\oplus\bigoplus_{j=1}^{h}\mathbb{C}\cdot E_{\kappa+1/2}(z,\chi_j).
\]
\textit{Proof.}
By Lemma 6.2, 11.1 and 11.2, it suffices to show that the constant terms in the Fourier expansion of $E'_{\kappa+1/2}(g,\chi_j)$ are linearly independent over $\mathbb{C}$ as functions in $g\in\widetilde{SL_2(\mathbb{A})}.$
We consider their values at $g=\mathbf{m}(a_1),...,\mathbf{m}(a_h).$
As we have seen in the last two sections, the constant term of $E'_{\kappa+1/2}(\mathbf{m}(a_i),\chi_j)$ is $f_{\chi_j}(\mathbf{m}(a_i)).$
Here $f_{\chi_j}$ is the function defined by equation (\ref{kgeq2}) in Section 8.
An easy calculation shows that 
\[
f_{\chi_j}(\mathbf{m}(a_i))=\chi_j(\mathfrak{a}_i)^{-1}f_{\chi_j}(1)\prod_{v<\infty}\frac{\alpha_{\psi_v}(\eta)}{\alpha_{\psi_v}(\eta a_{i,v})}|a_{i,v}|^{\kappa+1/2}
\]
where
\[
f_{\chi_j}(1)=2^{n(2\kappa-1/2)}|\mathfrak{D}|^{\kappa-1/2}\chi_j(\mathfrak{d})^{-1}\exp\left(\frac{(-1)^{\kappa+1}n\pi\sqrt{-1}}{4}\right)\neq0.
\]
Thus now we only have to show that the determinant of $(\chi_j(\mathfrak{a}_i)^{-1})_{0\leq i,j\leq h}$ is not zero.
But this is just the character table of the class group of $F.$\done\\


\section{An Example}
Let us consider the case $F=\mathbb{Q}(\sqrt{40})$ and $\kappa=2.$
The class number of $F$ is $h=2.$
Let the classes of ideals be $B_1$ and $B_2,$ where $B_1$ is the class consisting of principal ideals.
Also we let $\chi$ be the non-trivial character on the class group.\\
\hbox{}\quad By direct calculation, we have
\begin{align*}
60G_{5/2}(z,1)= &\\
1577
+70q
&+264q^2
+744q^{7\pm2\sqrt{10}}
+3850q^4
+3144q^5
+8640q^6
+\cdots, \\
60G_{5/2}(z,\chi)= &\\
1577
+24q
&+490q^2
+1750q^{7\pm2\sqrt{10}}
+2184q^4
+8470q^5
+8160q^6
+\cdots.
\end{align*}
\hbox{}\quad According to \cite{IkedaHiraga}, we have
\[
\dim_\mathbb{C}S^+_{\kappa+1/2}(\Gamma)=\dim_\mathbb{C}\mathbf{A}^\mathrm{CUSP}_{2\kappa}(\PGL_2(F)\backslash\PGL_2(\mathbb{A}_F)/\mathcal{K}_0)
\]
where $\mathbf{A}^\mathrm{CUSP}_{2\kappa}$ denotes the space of cuspidal automorphic forms of weight $2\kappa$ and $\mathcal{K}_0=\prod_{v<\infty}\PGL(\mathfrak{o}_v).$
In this case, it is known that the right hand side of this equation is 12, so by Corollary 11.3, we get
\[
\dim_\mathbb{C}M^+_{5/2}(\Gamma)=12+2=14.
\]
\hbox{}\quad One can check that the two Eisenstein series above are linear combinations of some already-known modular forms in $M^+_{5/2}(\Gamma).$
In fact, if we define the Eisenstein series of weight 2 by
\[
E_{2,i}(z)=\zeta_{F,i}(-1)/4+\sum_{\begin{smallmatrix}\xi\in\mathfrak{o}\\ \xi\succ0\end{smallmatrix}}\sigma_{1,i}(\xi)q^\xi
\]
for $i=1, 2,$ where
\[
\zeta_{F,i}(s)=\sum_{\begin{smallmatrix}\mathrm{integral }\end{smallmatrix}\mathfrak{c}\in B_i}N_{F/\mathbb{Q}}(\mathfrak{c})^{-s},\quad Re(s)>1,
\]
is being continued analytically to the whole complex plane and
\[
\sigma_{k,i}(\xi)=\sum_{\begin{smallmatrix}\mathrm{integral}\,\mathfrak{c}\in B_i\\ \mathfrak{c}\,|\,(\xi)\end{smallmatrix}}N_{F/\mathbb{Q}}(\mathfrak{c})^k,
\]
and two theta functions by
\[
\theta_1(z)=\sum_{\xi\in\mathfrak{o}}q^{\xi^2},\quad\theta_2(z)=\sum_{\xi\in\mathfrak{p}_5}q^{\xi^2/5},
\]
where $\mathfrak{p}_5^2=(5),$ then $E_{2,i}(4z)\theta_j(z)\in M^+_{5/2}(\Gamma)$ for $i, j=1, 2.$
Let 
\[
\begin{aligned}
&f_1(z)=E_{2,1}(4z)\theta_1(z),&f_2(z)=E_{2,1}(4z)\theta_2(z),\\
&f_3(z)=E_{2,2}(4z)\theta_1(z),&f_4(z)=E_{2,2}(4z)\theta_2(z).
\end{aligned}
\]
We have the identities
\[
\begin{aligned}
&172845227913G_{5/2}(z,1)=\frac{370247733672}{13}f_1-\frac{7861698464301}{91}f_3\\
&+16454261996f_4-\frac{6750047621}{26}T^+_{5/2}(\alpha_1^2)f_1+\frac{8395141929}{26}T^+_{5/2}(\alpha_1^2)f_2\\
&+\frac{37223824769}{104}T^+_{5/2}(\alpha_1^2)f_3-\frac{3375940624}{13}T^+_{5/2}(\alpha_1^2)f_4-\frac{649641221}{26}T^+_{5/2}(\alpha_2^2)f_1\\
&+\frac{1180397267}{26}T^+_{5/2}(\alpha_2^2)f_2+\frac{4022282847}{56}T^+_{5/2}(\alpha_2^2)f_3,\\
\end{aligned}
\]
\[
\begin{aligned}
&172845227913G_{5/2}(z,\chi)=-\frac{175639298994}{13}f_1+\frac{279576612332}{91}f_3\\
&+8260703363f_4-\frac{17803418247}{104}T^+_{5/2}(\alpha_1^2)f_1+\frac{24155608897}{104}T^+_{5/2}(\alpha_1^2)f_2\\
&+\frac{22793580805}{104}T^+_{5/2}(\alpha_1^2)f_3-\frac{7459199343}{52}T^+_{5/2}(\alpha_1^2)f_4+\frac{5253019763}{104}T^+_{5/2}(\alpha_2^2)f_1\\
&+\frac{4756228563}{104}T^+_{5/2}(\alpha_2^2)f_2-\frac{4528661307}{56}T^+_{5/2}(\alpha_2^2)f_3,\\
\end{aligned}
\]
where for any $f=\sum_{\xi}c(\xi)q^\xi\in M^+_{\kappa+1/2}(\Gamma)$ and $\alpha\in\mathfrak{o}$ such that $\mathfrak{p}_\alpha=(\alpha)$ is a prime ideal,
\[
\begin{aligned}
&T^+_{\kappa+1/2}(\alpha^2)f=\\
&\sum_{\begin{smallmatrix}\xi\equiv\square\,\mathrm{mod } 4\\ \xi\succ0\end{smallmatrix}}\left(c(\xi\alpha^2)+\left(\frac{(-1)^\kappa\xi}{\mathfrak{p}_\alpha}\right)N_{F/\mathbb{Q}}(\mathfrak{p}_\alpha)^{\kappa-1}c(\xi)+N_{F/\mathbb{Q}}(\mathfrak{p}_\alpha)^{2\kappa-1}c\left(\frac{\xi}{\alpha^2}\right)\right)q^\xi
\end{aligned}
\]
(we put $c(x)=0$ if $x\notin\mathfrak{o}$). Here $\alpha_1=3-2\sqrt{10}, \alpha_2=9-\sqrt{10}.$\\
\hbox{}\quad For more details of the Hecke operator $T^+_{\kappa+1/2}(\alpha^2)$ on $M^+_{\kappa+1/2}(\Gamma),$ one can consult \cite{Shimura:87}.
One can also directly see that $G_{\kappa+1/2}(z,\chi')$ is an eigenvector of $T^+_{\kappa+1/2}(\alpha^2)$ with eigenvalue $1+N_{F/\mathbb{Q}}(\mathfrak{p}_\alpha)^{2\kappa-1}.$

\end{document}